\documentstyle[amstex,12pt]{article}
\begin{document}
\renewcommand{\thefootnote}{\fnsymbol{footnote}}

\title{{\Large\bf GROUP RINGS
OF COUNTABLE NON-ABELIAN
LOCALLY FREE GROUPS ARE PRIMITIVE}}
%

\author{Tsunekazu Nishinaka \\
{\small Okayama Shoka University} 
{\small Okayama 700-8601} 
{\small Japan}}
\date{}
\footnotetext{it group ring; primitive; locally free group;
 ascending HNN-extension.\\
\quad\quad 2000 MSC:
6S34,\ 20C07, \ 20E25, \ 20E06
}
\maketitle
%

%
\baselineskip20pt
\newfont{\ef}{eufm10 scaled\magstep1}
\newfont{\efs}{eufm8}
\newfont{\bbl}{msbm10 scaled\magstep1}
\newfont{\cl}{cmsy10 scaled\magstep1}
\newfont{\cls}{cmsy8}

%
\newcommand{\St}{\mbox{\ef T}}
\newcommand{\XG}{\mbox{\ef X}}
\newcommand{\CG}{\mbox{\ef C}}
\newcommand{\cg}{\mbox{\ef c}}
\newcommand{\DG}{\mbox{\ef D}}
\newcommand{\UG}{\mbox{\ef U}}
\newcommand{\UGs}{\mbox{\efs U}}
\newcommand{\LG}{\mbox{\ef L}}
\newcommand{\MG}{\mbox{\ef M}}
\newcommand{\NGG}{\mbox{\ef N}}
\newcommand{\PG}{\mbox{\ef P}}
\newcommand{\WG}{\mbox{\ef W}}
\newcommand{\XGG}{\mbox{\ef X}}
\newcommand{\tl}{\mbox{\ef t}}
\newcommand{\hd}{\mbox{\ef h}}
\newcommand{\orn}{\mbox{\ef o}}
\newcommand{\ter}{\mbox{\ef t}}
\newcommand{\mx}{\mbox{\bf m}}
\newcommand{\XGGs}{\mbox{\efs X}}
\newcommand{\Zh}{\mbox{\bbl Z}}
\newcommand{\Nat}{\mbox{\bbl N}}
\newcommand{\Rd}{\mbox{\cl R}}
\newcommand{\Gh}{\mbox{\cl G}}
\newcommand{\Pw}{\mbox{\cl P}}
\newcommand{\Vt}{\mbox{\cl V}}
\newcommand{\Vts}{\mbox{\cls V}}
\newcommand{\Ed}{\mbox{\cl E}}
\newcommand{\Ghs}{\mbox{\cls G}}
\newcommand{\Sy}{\mbox{\cl R}}
\newcommand{\Sys}{\mbox{\cls R}}
\newcommand{\IG}{\mbox{\cl I}}
%
%
%
\newtheorem{theom}{{\bf Theorem}}[section]
\newtheorem{df}[theom]{{\bf Definition}}
\newtheorem{lemm}[theom]{{\bf Lemma}}
\newtheorem{coro}[theom]{{\bf Corollary}}
\newtheorem{propo}[theom]{{\bf Proposition}}
\newtheorem{rema}[theom]{{\bf Remark}}
\newtheorem{examp}[theom]{{\bf Example}}
\newtheorem{prb}[theom]{{\bf Problem}}

\vskip27pt

\begin{abstract}
We prove that every group ring of a non-abelian locally free group
which is the union of an ascending
sequence of free groups is primitive.
In particular,
every group ring of
a countable non-abelian locally free group is primitive.
In addition, by making use of the result,
we give a necessary and sufficient condition for group rings
of ascending HNN extensions of free groups
to be primitive,
which extends the main result in \cite{Ni07}
to the general cardinality case.
\end{abstract}
\vskip39pt

\section{INTRODUCTION}
\label{int}

A ring is (right) primitive
if it has a faithful irreducible (right) module.
Our purpose in this paper is to study
the primitivity of group rings of locally free groups.

A group is called locally free if
all of its finitely generated subgroups are free.
It is well known that there exist locally free groups
which are not free.
For example,
a properly ascending union of non-abelian free groups
of bounded finite rank is infinitely generated and Hopfian
(see \cite{Lyn}),
and so it is a locally free group which is not free,
where a group is Hopfian provided
that every surjective endomorphism of that group is an automorphism.
An example of uncountable non-free locally free groups
can be seen in Higman \cite{Hig}.
It has been seen
that a locally free group appears in a subgroup
of the fundamental group
of a three-dimensional manifold
(\cite{Freed}, \cite{Jam}, \cite{Kent}, \cite{Mask}).
The fundamental group of the mapping torus of the standard
2-complex of a free group $F$
with bounding maps the identity and
an injective endomorphism $\varphi$ of $F$,
which is called the ascending HNN extension of $F$
corresponding to $\varphi$,
also has a locally free group as a subgroup.

Recall that the ascending HNN extension $F_{\varphi}$ of $F$
corresponding to $\varphi$ has the presentation
$F_{\varphi}=\langle F, t| t^{-1}ft=\varphi(f)\rangle$.
The ascending HNN extension $F_{\varphi}$ of a free group $F$
is a well-studied class of groups.
For example, 
$F_{\varphi}$
is coherent (Feighn and Handel \cite{F-H}),
where a group is coherent
if its finitely generated subgroups are finitely presented.
If $F$ is finitely generated
then $F_{\varphi}$ is Hopfian
(Geoghegan, Mihalik, Sapir and Wise \cite{GMSW01}).
Moreover, 
Borisov and Sapir \cite{B-S}
have recently shown that it is residually finite.
The present author \cite{Ni07} has quite recently
shown that the group ring $KF_{\varphi}$ is semiprimitive
for any field $K$ and it is often primitive
provided that $F$ is a non-abelian countable free group.
In the proof, it was shown that
the group ring of a certain countable locally free group
is primitive, and therefore we posed a question:
Is it true that the group ring of any locally free group
is primitive?
In the present paper,
we shall give a partial answer to the question:

\begin{theom}\label{MainTh}
Let $G$ be a non-abelian locally free group
which has a free subgroup
whose cardinality is the same as that of $G$ itself.

$(1)$
Let $R$ be a domain (i.e. a ring with no zero divisors).
If $|R|\leq |G|$ then the group ring $RG$ is primitive.

$(2)$
If $K$ is a field then $KG$ is primitive.
\end{theom}

In particular,
every group ring of the union of an ascending
sequence of non-abelian free groups over a field
is primitive,
and so every group ring of
a countable non-abelian locally free group over a field
is primitive
(see Corollary \ref{MainCor}).
Theorem \ref{MainTh} corresponds to the result
obtained  by Formanek \cite{For}.
He showed that
every group ring $RG$ of a free product $G$ of non-trivial groups
(except $G=\Zh_{2}*\Zh_{2}$)
over a domain $R$ is primitive provided that
the cardinality of $R$
is not larger than that of $G$.
Motivated by the result,
the primitivity of some interesting rings and algebras
has been studied
(for example,
the primitivity of
free products of algebras by Lichtman \cite{Lic},
the primitivity of
group rings of amalgamated free products by Balogun \cite{Bal}
and the primitivity of semigroup algebras of free products
by Chaudhry, Crabb and McGregor \cite{C-C-M}).
In their papers
\cite{Bal} and \cite{Lic} (see also \cite{Law} and \cite{Ni07}),
the method established in \cite{For},
which is based on the construction of comaximal ideals,
has been applied to obtain the primitivity.
This Formanek's method is also available for our study.

In the present paper,
we state some graph-theoretic results
and apply them to Formanek's method.
For the sake of simplicity of the explanation,
we consider here  the group ring $KG$ over a field $K$.
For a non-zero element $u$ in $KG$,
let $\varepsilon(u)$ be an element in the ideal $KGuKG$
generated by $u$.
Then Formanek's method says that
$KG$ is primitive,
provided the right ideal $\rho$ generated by the elements $\varepsilon(u)+1$
for all non-zero $u$ in $KG$ is proper.
The main difficulty here is how to choose elements $\varepsilon(u)$'s
so as to make $\rho$ be proper.
That is, the chosen elements $\varepsilon(u)$'s must satisfy that
any finite sum of the form $r=\sum(\varepsilon(u)+1)v$
is not the identity element in $KG$, where $v$'s are elements in $KG$.
In general,
$\varepsilon(u)$'s and $v$'s are linear combinations
of elements of $G$;
say $f$'s and $g$'s are supports of $\varepsilon(u)$ and $v$
respectively.
If $r$ is not the identity element,
then it has at least one support of the form
$fg$ or $g$.
On the contrary, if $r$ is the identity element,
then almost all elements of the form $\alpha fg$
are vanished in $r$,
where $\alpha$ is a non-zero coefficient in $K$.
Then, what can we say
about supports $f$'s of $\varepsilon(u)$'s?
In order to consider this,
regarding the elements of the form
$fg$ or $g$ appeared in $r$ as vertices
and the equalities of their elements as edges,
we use a graph-theoretic method.
In section 3,
we define an R-graph and an R-cycle,
and show that under reasonable conditions,
two typical R-graphs have an R-cycle
(See Theorem \ref{PS3} and Theorem \ref{PS2}).
Roughly speaking,
if a suitable R-graph for vertices $fg$'s and $g$'s
has an R-cycle,
then it follows that
some supports $f_i$'s of $\varepsilon(u)$'s satisfy
the equation $f_1f_2^{-1}\cdots f_{n-1}f_n^{-1}=1$.
We should then choose $\varepsilon(u)$'s
so that their supports $f_i$'s never satisfy such equation.

Now,
if $\varphi$ is an automorphism, that is $\varphi(F)=F$,
then the ascending HNN extension
$F_{\varphi}$ is a cyclic extension of $F$.
On the other hand,
If $\varphi(F)\ne F$,
then it is a cyclic extension of a locally free group
(see Lemma \ref{fcc} (3)).
Therefore, by applying Theorem \ref{MainTh} to this case,
we can establish the primitivity of
the group ring $KF_{\varphi}$:

\begin{theom}\label{SubTh}
Let $F$ be a non-abelian free group.
Then the following are equivalent:

$(1)$ $KF_{\varphi}$ is primitive for a field $K$.

$(2)$ $|K|\leq |F|$
or
$F_{\varphi}$ is not virtually
the direct product $F\times \Zh$.

$(3)$ $|K|\leq |F|$ or $\triangle(G)=1$,
where $\triangle(G)$ is the FC center of $G$.

In particular,
if $F_{\varphi}$ is a strictly ascending HNN extension,
that is,
$\varphi(F)\ne F$,
then $KF_{\varphi}$ is primitive for any field $K$.
\end{theom}

This extends
the main result \cite[Theorem 1.1]{Ni07},
which was given for the countable case,
to the general cardinality case,
and follows the semiprimitivity of
$KF_{\varphi}$ with any cardinality
( see Corollary \ref{SubCor} ).


\section{PRELIMINARIES}
\label{pre}

Let $G$ be a group and $N$ a subgroup of $G$.
Throughout this paper,
we denote by $[G:N]$ the index of $N$ in $G$.
For a group $H$,
$G$ is said to be virtually $H$
if $H$ is isomorphic to $N$ and $[G:N]<\infty$.
If $g$ is an element of $G$,
we let $C_{N}(g)$ denote the centralizer of $g$ in $N$.
Let $C(G)$ be the center of $G$ and
$\triangle(G)$ the FC center of $G$,
that is $\triangle(G)=\{ g\in G\ |\ [G:C_{G}(g)]<\infty\}$.
Given a set $S$, let $|S|$ denote the cardinality of $S$.
If $S\subseteq G$ and $S=\{s_1,\cdots, s_m\}$,
$\langle S\rangle=\langle s_1,\cdots, s_m\rangle$
denotes the subgroup of $G$
generated by the elements of $S$.

The method of Formanek \cite{For}
based on the construction of comaximal ideals
plays also an important role in our study.
We shall give it as follows:
Let $G$ be an infinite group,
$R$ a ring with identity,
and $X$ a set with $|X|=|G|$.
Suppose that $|R|\leq |G|$.
Let $\psi$ be a bijection from $X$ to the elements of $RG$
except for the zero element.
For $x\in X$, let $\varepsilon^{*}(\psi(x))$
be a non-zero element in the ideal generated by $\psi(x)$ in $RG$,
and let $\rho=\sum_{x\in X}\varepsilon(\psi(x))RG$
be the right ideal of $RG$,
where $\varepsilon(\psi(x))=\varepsilon^{*}(\psi(x))+1$.

\begin{propo}\label{ForProp}\mbox{{\rm (See \cite{For})}}
If $\rho$ is proper then $RG$ is primitive.
\end{propo}

We shall use some basic results on free groups
in the last section.
For the details,
we refer the reader to
Lyndon and Schupp \cite{Lyn}.
The next assertions on locally free groups
are almost obvious.
For the sake of completeness,
we include a proof. 

\begin{lemm}\label{fg_b}
Let $G$ be a non-abelian locally free group.

$(1)$
The FC center of $G$ is trivial; thus $\Delta(G)=1$.

$(2)$
Let $S=\{v_{1},\cdots,v_{s}, w_{1},\cdots,w_{t}\}$
be a non-empty finite subset
of $G$ such that all elements in $S$ are non-trivial,
$v_{i}\ne v_{j}$ and $w_{i}\neq w_{j}$ if $i\ne j$.
Then for each $m>0$, there exist elements
$z_{1},\cdots,z_{m}\in G$
which satisfy

$(\mbox{{\rm i}})$
$v_{k}z_{l}w_{i}z_{l}
=v_{h}z_{n}w_{j}z_{n}$
if and only if $(k,l,i)=(h,n,j)$,

$(\mbox{{\rm ii}})$
for $p>0$,
if $\prod_{q=1}^{p}(z_{l_{q}}w_{i_{q}}z_{l_{q}})^{-1}
(z_{n_{q}}w_{j_{q}}z_{n_{q}})=1$,
then either $n_{q}= l_{(q+1)}$ for some $q\in\{1,\cdots,p-1\}$
or $(l_{q},i_{q})= (n_{q},j_{q})$ for some $q\in\{1,\cdots,p\}$.
\end{lemm}

\noindent
{\bf Proof.}\quad
Since $G$ is a non-abelian locally free group,
there exists a subset $X\subset G$ with $|X|>1$
such that $\langle X\rangle$ is 
freely generated by $X$ in $G$
and $\langle X\rangle\supseteq S$.
If $v\in S$ then $C_{\langle X\rangle}(v)$ is cyclic
(see \cite[Proposition 2.19]{Lyn}),
and so $[\langle X\rangle:C_{\langle X\rangle}(v)]$
is not finite,
which implies $v\notin\Delta(G)$.
Hence we see that $\Delta(G)=1$ and thus (1) holds.

Now, let $x_1,x_2\in X$ with $x_1\ne x_2$,
and let $n_{S}$ be the maximum length of
the words in $S$ on $X$,
where the length of a word $v$ is defined for
the reduced word equivalent to $v$  on $X$.
We set $n=2n_{S}$ and
$z_{l}=
x_1^{n+l}x_2x_1^{n+l},$
where $l=1,2,\cdots, m$.
Then it is easily verified
that the above $z_{1},\cdots,z_{m}$
satisfy (i) and (ii).
$\Box$
\vskip12pt

Let
$F$ be a non-abelian free group,
and $F_{\varphi}=\langle F, t| t^{-1}ft=\varphi(f)\rangle$
the ascending HNN extension of $F$ determined by $\varphi$.
It is easily verified that every element $g\in F_{\varphi}$
has a representation of the form $g=t^kft^{-l}$
where $k,l\geq 0$ and $f\in F$.
Combining this with some elementary observations on free groups,
we can show the following properties of $F_{\varphi}$:

\begin{lemm}\label{fcc}\mbox{{\rm (See \cite[Lemma 2.1, 2.2]{Ni07})}}
Let $F$ be a non-abelian free group.

$(1)$
$\Delta(F_{\varphi})=C(F_{\varphi})$.

$(2)$
The following are equivalent:

\indent
\mbox{{\rm (i)}}
$C(F_{\varphi})\ne 1$.

\indent
\mbox{{\rm (ii)}}
There exist $n>0$ and $f\in F$
such that $C(F_{\varphi})=\langle t^{n}f\rangle$.

\indent
\mbox{{\rm (iii)}}
$F_{\varphi}$ is virtually the direct product
$F\times \Zh$.

When this is the case,
$\varphi$ is an automorphism of $F$;
thus $\varphi(F)=F$.

$(3)$
For a non-negative integer $i$,
let $F_{i}$ be the subgroup of $F_{\varphi}$
generated by $\{ t^{i}ft^{-i}\ |\ f\in F\}$.
Then $F_1\subseteq F_2\subseteq\cdots \subseteq F_i\subseteq\cdots$
is an ascending chain of free groups,
and $F_{\infty}=\bigcup_{i=1}^{\infty}F_{i}$ is a normal
subgroup of $F_{\varphi}$.
\end{lemm}

The next two lemmas are basic results on group rings.
We refer the reader to Passman \cite{Ps77}
for a more detailed discussion of these questions.

\begin{lemm}\label{PBR}
Let $K$ be a field,
$G$ a group and $N$ a subgroup of $G$.

\noindent
$(1)$\mbox{{\rm (See \cite[Theorem 1]{Zal})}}
Suppose that $N$ is normal.
If $\triangle(G)$ is trivial and $\triangle(G/N)=G/N$,
then $KN$ is primitive implies $KG$ is primitive.

\noindent
$(2)$\mbox{{\rm (See \cite[Theorem 3]{Rosn})}}
If $\triangle(G)$ is torsion free abelian and $[G:N]$ is finite,
then $KN$ is primitive implies $KG$ is primitive.
\end{lemm}

\begin{lemm}\label{Ps}\mbox{{\rm (See \cite[Theorem 2]{Ps73})}}
Let $K'$ be a field and $G$ a group.
If $\triangle(G)$ is trivial and $K'G$ is primitive,
then for any field extension $K$ of $K'$,
$KG$ is primitive.
\end{lemm}

Formanek \cite{For}
asserts that if $G$ is the direct product of a free group $F$
and the infinite cyclic group $\langle t\rangle$,
then $KG$ is primitive
if and only if the cardinality of $K$ is not larger than that of $F$.
Combining this with Lemma \ref{fcc} and \ref{PBR} (2),
we have

\begin{lemm}\label{hnn}\mbox{{\rm (See \cite[Theorem 1.1 (i)]{Ni07})}}
Let $F$ be a non-abelian free group,
and suppose that $\varphi(F)=F$,
that is, $KF_{\varphi}$is the cyclic extension of $F$ by $\langle t\rangle$.
Then the following are equivalent:

$(1)$ $KF_{\varphi}$ is primitive for a field $K$.

$(2)$ $|K|\leq |F|$
or
$F_{\varphi}$ is not virtually
the direct product $F\times \Zh$.

$(3)$ $|K|\leq |F|$ or $\triangle(G)=1$.
\end{lemm}

%

\section{GRAPHICAL INVESTIGATION}

Let $KG$ be the group ring of a group $G$ over a field $K$,
and let $a=\sum_{i=1}^{m}\alpha_if_i$
and $b=\sum_{i=1}^{n}\beta_ig_i$ be in $KG$
$(\alpha_i\ne 0, \beta_i\ne 0)$.
If $ab=0$ then
for each $f_ig_j$, there exists $f_{p}g_{q}$
such that $f_ig_j=f_{p}g_{q}$.
Suppose that the following $k$ equations hold;
$f_1g_1=f_2g_2$, $f_3g_2=f_4g_3$, $\cdots,$
$f_{2k-3}g_{k-1}=f_{2k-2}g_{k}$
and $f_{2k-1}g_k=f_{2k}g_{1}$.
Then
we can regard
the above equations as forming a kind of cycle,
and they imply
$f_1^{-1}f_2\cdots f_{2k-1}^{-1}f_{2k}=1$.
That is, the above equations give us a information
on supports of $a$.
We can use this idea for a more general case;
$a_1b_1+\cdots+a_nb_n\in K$ for $a_i,b_i\in KG$
with $a_i=\sum\alpha_{ij}f_{ij}$ and $b_i=\sum\beta_{ik}g_{ik}$.
To do this,
regarding the elements $f_{ij}g_{ik}$ appeared in $a_ib_i$ as vertices
and the equalities of their elements as edges,
we use a graph-theoretic method.

In this section,
we shall define an R-graph and an R-cycle,
and show that under reasonable conditions,
two typical R-graphs have an R-cycle.
One is called an R-colouring R-graph
and another is called an R-simple R-graph
which is a special case of an R-colouring.
These two results, Theorem \ref{PS3} and Theorem \ref{PS2},
are used to prove our main theorem in the next section.

Throughout this section,
$\Gh=(V,E)$ denotes a simple graph;
a finite undirected graph
which has no multiple edges or loops,
where $V$ is the set of vertices
and $E$ is the set of edges.
For terminology and notations not defined here, 
we refer to Bondy and Murty \cite{B-M}.
A finite sequence $v_0e_1v_1\cdots e_pv_p$
whose terms are alternately
elements $e_q$'s in $E$ and $v_q$'s in $V$
is called a path of length $p$ in $\Gh$
if $v_{q-1}v_q=e_q\in E$
and $v_q\ne v_{q'}$
for any $q,q'\in\{0,1,\cdots,p\}$ with $q\ne q'$;
simply denoted by $v_0v_1\cdots v_p$.
Two vertices $v$ and $w$ of $\Gh$ are said to be connected
if there exists a path from $v$ to $w$ in $\Gh$.
Connection is an equivalence relation on $V$,
and so there exists a decomposition of $V$ into subsets $C_i$'s
$(1\leq i\leq m)$ for some $m>0$ such that $v, w\in V$
are connected if and only if both $v$ and $w$ belong to
the same set $C_i$;
each $C_i$ is called a (connected) component of $\Gh$.
Any graph is a disjoint union of components.



\begin{df}\label{Rgrh}
Let $\Gh=(V,E)$ and $\Gh^{*}=(V,E^{*})$ be simple graphs
with the same vertex set $V$.
For $v\in V$,
let $U(v)$ be the set consisting of all neighbours
of $v$ in $\Gh^{*}$ and $v$ itself:
$U(v)=\{w\in V\ |\ vw\in E^{*}\}\cup\{v\}$.
A triple $(V,E, E^{*})$ is an R-graph
(for a relay-like graph)
if it satisfies the following condition $(\mbox{{\rm R}})$:
\vskip6pt

$(\mbox{{\rm R}})$
If $v\in V$
and $C$ is a component of $\Gh$,
then $|U(v)\cap C|\leq 1$.
\vskip6pt

\noindent
That is,
each $U(v)$ has at most one vertex from each component of $\Gh$.
If $\Gh$ has no isolated vertices,
that is, 
if $v\in V$ then
$vw\in E$ for some $w\in V$,
then R-graph $(V,E, E^{*})$ is called a proper R-graph.
\end{df}

We call $U(v)$ the R-neighbour set of $v\in V$,
and set $\UG=\{ U(v)\ |\ v\in V\}$.
For $v,w\in V$ with $v\ne w$,
it may happen that $U(v)=U(w)$,
and so $|\UG|\leq |V|$ generally.
If $w, w'\in U(v)$ then the minimum length of paths from $w$ to $w'$
in $\Gh^{*}$ is at most $2$.
Moreover, $|U(v)\cap C|>1$ for some $v\in V$
and for some component $C$ of $\Gh$
if and only if there exists a path from $w$ to $w'$ for some $w,w'\in U(v)$
in $\Gh$.
Hence we have

\begin{propo}\label{REC}
In the definition \ref{Rgrh},
the condition \mbox{{\rm (R)}} is equivalent to
each of the following conditions:

$(\mbox{{\rm R}}')$
If $C$ is a component of $\Gh$ and $C'$ is a component of $\Gh^{*}$
and if there exist $v,w\in C\cap C'$ with $v\ne w$,
then the length of any path from $v$ to $w$ in $\Gh^{*}$ is longer than $2$.

$(\mbox{{\rm R}}'')$
If $U\in \UG$ and $v,w\in U$,
then there exist no paths from $v$ to $w$;
if $v_0v_1\cdots v_m$ is a path in $\Gh$,
then $\{v_0,v_m\}\not\subseteq U$ for any $U\in\UG$.
\end{propo}

By $(\mbox{{\rm R}}')$,
in particular,
if $vw\in E^{*}$
then $vw\not\in E$.
We say $\Gh=(V,E)$ to be the base graph of $\Sy=(V,E,E^{*})$.
If $E=\emptyset$
then $\Sy$
is called the empty graph; denoted by $\Sy=\emptyset$.
Clearly, if $\Sy$ is non-empty
then $|\UG|>1$. 

In what follows,
let $\Sy=(V, E, E^{*})$ be a non-empty R-graph
with $\Gh=(V,E)$ and $\Gh^{*}=(V,E^{*})$.
For $W\subseteq V$,
we define $E_{W}$ and $E_{W}^{*}$ by
$$
\begin{array}{lll}
E_{W}&=
\{ww'\ |\ w,w'\in W \mbox{ and either } ww'\in E\\
&\quad\ \mbox{ or }\
wv_1\cdots v_mw' \mbox{ is a path in } \Gh\ \mbox{ for some } 
v_i\in V\setminus W\},\\
E_{W}^{*}&=
\{ww'\ |\ w,w'\in W \mbox{ and } ww'\in E^{*}\}.\\
\end{array}$$
$E_{W}^{*}$ is simply the edges of
the subgraph of $\Gh^{*}$ generated by $W$.
$E_{W}$ means that if $C$ is a component of $\Gh$
with $C\cap W\ne \emptyset$,
then $C\cap W$ also becomes a component of $\Gh_W=(W,E_{W})$.
It is obvious that
$\Sy_{W}=(W, E_W, E_{W}^{*})$
is an R-graph.
We call $\Sy_{W}$
the R-subgraph of $\Sy$
generated by $W$,
and set $\UG_{W}=\{ U_{W}(v)=U(v)\cap W\ |\ v\in W,\ U(v)\in\UG\}$.
If $E_{W}$ coincides with
$\{ww' \mid w,w'\in W, ww'\in E\}$,
then $\Sy_{W}$ is simply called the subgraph of $\Sy$
generated by $W$.

The degree $d_W(v)$ of $v\in W$ in $\Sy_W$
is the number of edges of $\Gh_W$ incident with $v$;
$d_V(v)$ is simply denoted by $d(v)$.
For $X\subseteq W\subseteq V$, $I_W(X)$ denotes
$\{x\in X\ |\ d_W(x)=0\}$;
$I_V(W)$ is simply denoted by $I(W)$.
In general, $|I_W(W)|\geq |I(W)|$.
In fact, we have

\begin{lemm}\label{SimSup}
Let $W\subseteq V$ and $W^{c}=V\setminus W$.
Then
$$0\leq |I_W(W)|-|I(W)|\leq |W^{c}|-|I(W^{c})|.$$
The right side equality holds if and only if
for each $v\in W^{c}\setminus I(W^{c})$, $d(v)=1$ and
there exists $w\in W$ with $d(w)=1$
such that $vw\in E$.
\end{lemm}

\noindent
{\bf Proof.}\quad
Let $X=I_W(W)\setminus I(W)$ and $Y=W^{c}\setminus I(W^{c})$.
Since $0\leq |X|$ is obvious,
it suffices to show that $|X|\leq |Y|$ and the assertion on equality in the statement
is true.

Note that $w\in X$ if and only if $d_W(w)=0$ and $d(w)\ne 0$.
Therefore, if $w\in X$ then there exists $v\in Y$ such that $vw\in E$.
Then $vw'\not\in E$ for any $w'\in X$ with $w'\ne w$.
In fact,
if $vw'\in E$ then $ww'\in E_W$,
which implies  a contradiction that $w,w'\not\in I_{W}(W)$.
Hence $|X|\leq |Y|$.
Since $|X|<\infty$ and $|Y|<\infty$,
if $|X|=|Y|$, then for each $w\in X$
( resp. for each $v\in Y$ ) there exists only one $v\in Y$
( resp. $w\in X$ ) such that $vw\in E$.
From this, if for $w\in X$, $d(w)>1$
then $ww'\in E$ for some $w'\in W$ with $w\ne w'$,
but this implies
$w,w'\not\in I_W(W)$, a contradiction.
Hence $d(w)=1$ for all $w\in X$.
On the other hand,
if $d(v)>1$ for some $v\in Y$
then $vv'\in E$ for some $v'\in Y$ with $v\ne v'$.
For these $v,v'\in Y$,
as mentioned above,
there exists $w,w'\in X$ with $w\ne w'$
such that $vw,v'w'\in E$,
which implies a contradiction $w,w'\not\in I_{W}(W)$,
because $wvv'w'$ is a path in $\Sy$ and so $ww'\in E_W$.
We have therefore that $d(v)=1$ for all $v\in Y$.

The converse assertion on equality is obvious.
$\Box$
\vskip12pt

In what follows,
for $\pi=v_0v_1\cdots v_p$ a path in $\Gh$,
the origin $v_0$ of $\pi$ and the terminus $v_p$ of $\pi$
are denoted by $\orn(\pi)$ and $\ter(\pi)$
respectively. 

\begin{df}\label{Rcyc}
Let $p>1$ and let $\pi_q$ be a path in $\Gh$
with $v_q=\orn(\pi_q)$ and $w_q=\ter(\pi_q)$ $(1\leq q\leq p)$.
Then a sequence $(\pi_1, \pi_2, \cdots, \pi_p)$
is an R-path of length $p$ in $\Sy$
if it satisfies the following conditions
\mbox{\rm (i)} and \mbox{\rm (ii)}:

\mbox{\rm (i)}
All of $v_q$'s and $w_q$'s are different from each other,

\mbox{\rm (ii)}
$w_{q}v_{q+1}\in E^{*}$ for $1\leq q\leq p-1$.

\noindent
If, in addition, 
it satisfies the following condition \mbox{\rm (iii)},
then it is an R-cycle of length $p$ in $\Sy$:

\mbox{\rm (iii)}
$w_{p}v_{1}\in E^{*}$.

In particular,
if the length of $\pi_q$ is $1$,
that is, $\pi_q=e_q\in E$ for all $q$,
then $(e_1, e_2, \cdots, e_p)$ is called an R-cycle consisting of edges.
\end{df}

It is obvious that $\Sy$ has an R-cycle
if the R-subgraph $\Sy_{W}$
has an R-cycle for some $\emptyset\ne W\subseteq V$.
A proper R-graph
$\Sy$ is called a clique R-graph,
provided that the base graph $\Gh=(V,E)$ is a clique graph; thus
$uv,vw\in E$ implies $uw\in E$.
Note that $\Sy$ is a clique R-graph
if and only if every component of $\Gh$
is a complete graph.
Hence,
if a clique R-graph has an R-cycle then
it also has an R-cycle consisting of edges.

In what follows,
$\NGG(\Sy)=\{U\in \UG\ |\ |U|=1\}$.
We note, if a sequence
$(\pi_1,\cdots, \pi_p)$
is an R-cycle in $\Sy$,
then neither $U(\orn(\pi_q))$ nor $U(\ter(\pi_q))$
is in $\NGG(\Sy)$
for all $1\leq q\leq p$.

We consider the following condition (UC) for $\UG$ in $\Sy$:
\vskip12pt

\noindent
\mbox{{\rm (UC)}}\quad
$\mbox{For each } U \mbox{ and }
U' \mbox{ in } \UG,
\mbox{ either } U\cap U'=\emptyset
\mbox{ or } |U\cap U'|>1.$
\vskip12pt

If $\UG$ satisfies (UC)
then $\NGG(\Sy)=\emptyset$,
because $|U|=|U\cap U|>1$ for each $U\in\UG$.

\begin{lemm}\label{PS1}
Let $\Sy$ be a proper R-graph.
If $\UG$ satisfies \mbox{{\rm (UC)}}
then $\Sy$ has an R-cycle.
\end{lemm}

\noindent
{\bf Proof.}\quad
Let $v_1\in V$.
Since $\Sy$ is proper,
there exists $w_1\in V$ such that $e_1=v_1w_1\in E$.
Since $|U(w_1)|>1$,
there exists $v_2\in V$ such that $w_1v_2\in E^{*}$.
Then $v_2v_1\not\in E^{*}$ by (R).
Since $\Sy$ is proper again,
there exists $w_2\in V$
such that $e_2=v_2w_2\in E$,
where $w_2\ne w_1$ and $w_2\ne v_1$
because of (R).
We see then that $(e_1, e_2)$ satisfies both of (i) and (ii)
in Definition \ref{Rcyc}; thus it is an R-path in $\Sy$.
If $U(w_2)\not\subseteq \{v_1, w_1\}$,
then we can proceed with this procedure.
Since 
$\Sy$ is a finite graph, this procedure terminates
in a finite number of steps,
say, exactly $p$ steps; that is,
there exist $p>1$
and a sequence $\sigma=(e_1,\cdots ,e_{p-1})$ of edges
with $e_q=v_qw_q$ such that
$\sigma$ is an R-path in $\Sy$,
and in addition, there exists $e_p=v_pw_p\in E$ such that $w_{p-1}v_p\in E^{*}$
and
$$U(w_p)\subseteq \{v_q,w_q\ |\ 1\leq q\leq p-1\}.$$
If there exists $q\in \{1\leq q\leq p-1\}$
such that $v_q\in U(w_p)$,
then the sequence $\sigma$ contains an R-cycle.
In fact, if $w_p=v_q$ for some $q\in\{1,\cdots,p-1\}$
then $q< p-1$ by (R),
and then $(\pi_q,e_{q+1},\cdots ,e_{p-1})$ is an R-cycle,
where $\pi_q=v_pv_qw_q$.
If $v_q\in U(w_p)\setminus \{w_p\}$ for some $q\in\{1,\cdots,p-1\}$;
thus $w_pv_q\in E^{*}$,
then $(e_{q},\cdots,e_p)$
is an R-cycle.
Therefore,
we may assume that for each $q\in \{1\leq q\leq p-1\}$,
$v_q\not\in U(w_p)$, that is
$$w_pv_q\not\in E^{*}\quad (1\leq q\leq p-1).$$
Then $w_q\in U(w_p)$ for some $q\in \{1\leq q\leq p-1\}$.
If $w_p=w_q$,
$(e_{q+1},\cdots,e_p)$ is an R-cycle because $w_pv_{q+1}=w_qv_{q+1}\in E^{*}$.
We may assume therefore that $w_pw_q\in E^{*}$ and
$q$ is minimal with this property;
thus $q=min\{1\leq q'\leq p-1\ |\ w_pw_{q'}\in E^{*}\}$.
Note that $q<p-1$ by (R).
Since $U(w_p)\cap U(v_{q+1})\supseteq \{w_q\}\ne\emptyset$,
we have that $|U(w_p)\cap U(v_{q+1})|>1$ by the hypothesis (UC).
Since $w_pv_{q+1}\not\in E^{*}$,
there exists $q'\in \{1\leq q\leq p-2\}$ with $q'\ne q$
such that $w_pw_{q'}\in E^{*}$ and $v_{q+1}w_{q'}\in E^{*}$.
By the minimality of $q$, $q<q'$, in fact, $q+1<q'$ by (R).
Since $w_{q'}v_{q+1}\in E^{*}$,
we see then that $(e_{q+1},\cdots,e_{q'})$ is an R-cycle.
$\Box$
\vskip12pt

In the above lemma,
we cannot replace the condition (UC)
by $\NGG(\Sy)=\emptyset$.
For instance, we have the following example:

\begin{examp}\label{ex1}
Let  $\Gh=(V,E)$ be the base graph
with $V=\{v_1,\cdots,v_{10}\}$
and $E=\{ v_1v_3,$ $v_2v_5, v_4v_7, v_6v_9, v_8v_{10}\}$.
Let $E^{*}=\{v_1v_2, v_3v_4, v_4v_5, v_6v_7, v_7v_8, v_9v_{10}\}$.
Then we have that $U(v_1)=U(v_2)=\{v_1,v_2\}$,
$U(v_3)=\{v_3,v_4\}$, $U(v_5)=\{v_4,v_5\}$,
$U(v_4)=\{v_3,v_4,v_5\}$,
$U(v_6)=\{v_6,v_7\}$, $U(v_8)=\{v_7,v_8\}$,
$U(v_7)=\{v_6,v_7,v_8\}$
and $U(v_9)=U(v_{10})=\{v_9,v_{10}\}$.
In this case, $\Sy=(V,E,E^{*})$ is a non-empty R-graph
and $\NGG(\Sy)=\emptyset$
but it has no R-cycles.
\end{examp}

Let $\CG_n(V)=
\{ V_1,\cdots,V_n\}$ be the set of components
of $\Gh^{*}=(V,E^{*})$.
For $1\leq i\leq n$ and for $v,v'\in V_i$,
define $v\simeq v'$ by $U^{o}(v)=U^{o}(v')$,
where $U^{o}(v)=U(v)\setminus \{v \}$;
thus $U^{o}(v)$ is the set of neighbours of $v$ in $\Gh^{*}$.
Clearly, $\simeq$ is an equivalence relation on $V_i$.
If $\cg(V_i)=\{ V_{i1},\cdots,V_{il_{i}}\}$
is the set of equivalence class of $V_i$,
then $V_i$ is the disjoint union of non-empty $V_{ij}$'s.
We can easily see that
\begin{equation}
\begin{array}{cc}
\mbox{if $v\in V_i$, there exists } S\subseteq \cg(V_i)\setminus \{V_{ij}\}
\mbox{ such that }\\
U^{o}(v)=\bigcup_{W\in S}W
\mbox{ for all } v\in V_{ij}.\\
\end{array}
\end{equation}
In particular, for each $v,w\in V_{ij}$,
$w\not\in U^{o}(v)$.
If for each $i\in \{1,\cdots,n\}$ and for each $v,v'\in V_i$,
$U(v)\cap U(v')\ne\emptyset$,
then $\CG_n(V)$ is called a colouring of $\Sy$.
Note that $\CG_n(V)$ is a colouring
if and only if for each $i$ and for each $v, v'\in V_i$,
there exists $U\in\UG$ such that $v,v'\in U$,
and so if $\CG_n(V)$ is a colouring
then for each $i$ and for each $v,v'\in V_i$,
$vv'\not\in E$ by (R$''$);
thus, in this case,
$\CG_n(V)$ is a colouring of the base graph $\Gh$.

We here consider the following condition (RC) for $\cg(V_i)$
which is stronger than \mbox{{\rm (1)}} above:
\vskip12pt

\noindent
\mbox{{\rm (RC)}}\quad
For each $v\in V_{ij}$,
$U^{o}(v)=V_{i}\setminus V_{ij}$.
\vskip12pt

\noindent
That is, $\cg(V_i)$ satisfies (RC)
if and only if $\Gh^{*}(V_i, E_{V_i}^{*})$ is a complete $k$-partite graph
$K_{l_1,\cdots,l_k}$,
where $k=|\cg(V_i)|$ and $l_j=|V_{ij}|$.

Let consider the graph $G_i=(\cg(V_i),\Ed)$
with the vertex set $\cg(V_i)=\{ V_{i1},\cdots,V_{il_{i}}\}$
and the edge set $\Ed=\{V_{ij}V_{ik}\ |\ j\ne k, \mbox{ for }
v\in V_{ij}, U^{o}(v)\supseteq V_{ik}\}$.
If $U^{o}(v)\supseteq V_{ik}$ for $v\in V_{ij}$
then $U^{o}(v')\supseteq V_{ij}$ for $v'\in V_{ik}$.
Combining this with (1),
we see that the above definition of $\Ed$ is well defined.
Then, the neighbour set of $V_{ij}$ in $G_i$
is $\cg(V_i)\setminus \{V_{ij}\}$
if and only if for each $v\in V_{ij}$,
$U^{o}(v)=V_{i}\setminus V_{ij}$ in $\Sy$,
and therefore,
$\cg(V_i)$ satisfies (RC)
if and only if $G_i$ is isomorphic to
the complete graph $K_{l_i}$ on $l_i$ vertices.
Now, the definition of $\CG_n(V)$
implies that $G_i$ is a connected graph,
and the definition of $\cg(V_i)$ implies
that for each distinct vertices $V_{ij},V_{ik}$ in $\cg(V_i)$,
the neighbour set of $V_{ij}$ does not coincide with that of $V_{ik}$.
It easily follows from these above that
$G_i$ is isomorphic to the complete graph $K_{l_i}$,
provided $l_i\leq 3$.
Hence we have

\begin{rema}\label{rem1}
If $|\cg(V_i)|\leq 3$,
then $\Gh^{*}_{V_i}=(V_i, E_{V_i}^{*})$ is a complete $k$-partite graph;
thus
$\cg(V_i)$ satisfies \mbox{{\rm (RC)}}.
\end{rema}

In Example \ref{ex1},
we can set that
$V_{11}=\{v_1\}$, $V_{12}=\{v_2\}$,
$V_{21}=\{v_3,v_5\}$, $V_{22}=\{v_4\}$,
$V_{31}=\{v_6,v_8\}$, $V_{32}=\{v_7\}$,
$V_{41}=\{v_9\}$, $V_{42}=\{v_{10}\}$
and $V_i=\bigcup_{j}V_{ij}$.
Then $\Gh^{*}_{V_1}\simeq \Gh^{*}_{V_4}\simeq K_{1,1}$
and $\Gh^{*}_{V_2}\simeq \Gh^{*}_{V_3}\simeq K_{1,2}$.

Let $l_i=|\cg(V_i)|$.
If $l_i>3$,
$G_i$ need not be isomorphic to the complete graph $K_{l_i}$
and thus
$\cg(V_i)$ need not satisfy (RC) (See Example \ref{ex2} below).
Our purpose is to make use of results on R-graphs for
proving our main theorem in the next section.
From this point of view,
it suffices to consider the case when $l_i=3$.
However, in the following consideration,
we only need the assumption (RC);
we need not to assume the condition $l_i\leq 3$.
We therefore define $\CG_n(V)$ to be an R-colouring of $\Sy$
if for each $i$, $\cg(V_i)$ satisfies (RC),
and investigate when R-colouring R-graphs have an R-cycle.
As a result,
we can use R-graph theory to analyze more general case
than the one of our main theorem (See Corollary \ref{subMainCor}).

\begin{df}\label{Rcol}
Let $\Sy$ be an R-graph with $\CG_n(V)=\{V_1,\cdots,V_n\}$.
If for each $1\leq i\leq n$,
$\Gh^{*}_{V_i}=(V_i, E_{V_i}^{*})$ is a complete $k$-partite graph;
thus
$\cg(V_i)$ satisfies \mbox{{\rm (RC)}},
then $\CG_n(V)$ is an R-colouring of $\Sy$
or $\Sy$ is an R-colouring R-graph with $\CG_n(V)$.
\end{df}


If $\CG_n(V)$ is an R-colouring of $\Sy$,
then it is a colouring of $\Sy$.
As has been mentioned in Remark \ref{rem1},
in case of $l_i=|\cg(V_i)|\leq 3$,
an R-graph is always an R-colouring.
However,
in case of $l_i> 3$,
it is not true. In fact,
a colouring need not be an R-colouring.
If $\CG_n(V)$ is a colouring of $\Sy$,
then each $V_{ij},V_{ik}\in\cg(V_i)$,
there exists a path $\pi$ in $G_i=(\cg(V_i),\Ed)$
whose length is $1$ or $2$
such that the origin $\orn(\pi)=V_{ij}$
and the terminus $\ter(\pi)=V_{ik}$.
Hence, for example, if $l_i=4$,
$G_i$ is isomorphic to either
the complete graph $K_4$ or the graph described in
the following example:

\begin{examp}\label{ex2}
Let $V_i=\{v_1,v_2,v_3,v_4\}$,
$U(v_1)=\{v_1,v_2\}$,
$U(v_2)=V_i$,
$U(v_3)=U(v_4)=\{v_2, v_3, v_4\}$.
In this case, we have that
$V_{ij}=\{v_j\}$ for $1\leq j\leq 4$;
thus $|\cg(V_i)|=4$.
Then,
$v_2\in U(v_j)$ for all $1\leq j\leq 4$; thus
$U(v_p)\cap U(v_q)\ne\emptyset$ for each $p,q\in\{1,2,3,4\}$,
but for $v_1\in V_{i1}$,
$U^{o}(v_1)=\{v_2\}=V_{i2}\ne V_i\setminus V_{i1}$.
Hence, $\cg(V_i)$ satisfies the colouring condition
but it fails to satisfy \mbox{{\rm (RC)}}, and certainly,
in the graph $G_i=(\cg(V_i),\Ed)$,
we see that
$\Ed$ coincides with
$\{V_{i1}V_{i2},V_{i2}V_{i3},V_{i3}V_{i4},V_{i4}V_{i2}\}$;
thus $G_i$ is not isomorphic to the complete graph $K_4$.
\end{examp}

Let $\Sy$ be an R-colouring R-graph with $\CG_n(V)$,
where $\CG_n(V)=\{ V_1,\cdots,V_n\}$
and $\cg(V_i)=\{ V_{i1},\cdots,V_{il_{i}}\}$.
For $W\subseteq V_i$,
We denote by $\mx(W)$
the maximum number in $\{|W\cap V_{ij}| \ |\ 1\leq j\leq l_i\}$
and by $J_{W}$ the set $\{j\ |\ W\cap V_{ij}\ne\emptyset\}$.
In general, $\mx(V_i)\geq \mx(W)$,
and clearly, if $\mx(W)=1$ then $\Gh^{*}_{W}=(W,E^{*}_{W})$ is a complete graph,
that is,
$U_{W}(v)=W$ for all $v\in W$,
and also if  $|J_{W}|>2$ then $|U_{W}(v)\cap U_{W}(w)|>1$
for all $v,w\in W$.
Suppose $|W|>\mx(W)+1$.
Then $l_{i}>1$ 
and $|J_{W}|>1$.
In this case,
if $|J_{W}|=2$, say $J_{W}=\{1,2\}$, and $\mx(W)=|W\cap V_{i1}|$,
then $\mx(W)>1$ and $|W\cap V_{i2}|=|W|-\mx(W)>1$,
which implies that $|U_{W}(v)\cap U_{W}(w)|>1$
for all $v,w\in W$.
Hence we have

\begin{rema}\label{rc1}
Let $W\subseteq V_i$.

\mbox{{\rm (i)}}
If $\mx(W)=1$ then $U_{W}(v)=W$ for all $v\in W$.

\mbox{{\rm (ii)}}\
If $|J_{W}|>2$ then $|U_{W}(v)\cap U_{W}(w)|>1$ for all $v,w\in W$.

\mbox{{\rm (iii)}}\
If $|W|>\mx(W)+1$
then $|U_{W}(v)\cap U_{W}(w)|>1$ for all $v,w\in W$.
\end{rema}

Recall that for $W\subseteq V$, $I(W)$ denotes
$\{w\in W\ |\ d_V(w)=0\}$.

\begin{theom}\label{PS3}
Let $n>1$, and let
$\Sy=(V,E,E^{*})$ be an R-colouring R-graph with
$\CG_n(V)=\{ V_1,\cdots,V_n\}$.
Suppose that $|V_i|\geq 2\mx(V_i)+1$ for each $i\in\{1,\cdots, n\}$.
If $|I(V)|\leq n$
then $\Sy$ has an R-cycle.
\end{theom}

\noindent
{\bf Proof.}\quad
Let $W=V\setminus I(V)$,
$W_i=V_i\setminus I(V_i)$
and $m_i=\mx(V_i)$ for $i=1,\cdots,n$.

We prove the statement above
by induction on $n$.
First, let $n=2$.
By the assumption,
$$|W_i|=|V_i|-|I(V_i)|\geq 2m_i+1-2\geq m_i>0\quad (i=1,2),$$
and so the R-subgraph $\Sy_{W}$
( in this case, it is simply the subgraph generated by $W$ )
is non-empty; thus $|J_{W_i}|>0$ for $i=1,2$.
Moreover,
\begin{equation}
|W|=|V|-|I(V)|\geq |V|-2\geq \sum_{i=1}^{2}(2m_i+1)-2=2(m_1+m_2).
\end{equation}
If $|J_{W_i}|=1$ for $i=1,2$,
then $|W|\leq m_1+m_2$, 
which contradicts  (2) above,
and so $|J_{W_i}|>1$ for $i=1$ or $i=2$.
If $|J_{W_1}|=1$ then
$|J_{W_2}|>1$ and $|W_1|<|W_2|$ 
because $|W_1|\leq m_1< m_1+2m_2\leq |W_2|$ by (2).
Since $\Sy_{W}$ is proper and $|W_1|<|W_2|$,
there exist $v\in W_1$ and $v_1,v_2\in W_2$ with $v_1\ne v_2$
such that $vv_1, vv_2\in E$.
However, since $\Sy$ is an R-colouring R-graph and $|J_{W_2}|>1$,
there exists $U\in \UG_{W_2}$ such that $v_1,v_2\in U$;
this contradicts (R).
We see therefore that $|J_{W_i}|>1$
for both $i=1$ and $i=2$,
and also that
$d(v)=1$ for all $v\in W$
and $|W_1|=|W_2|$.
Again by (2), we have that $|W_i|\geq m_1+m_2$.
In particular, $|W_i|\geq 2$ $(i=1,2)$.
If $m_1=1$ or $m_2=1$,
say $m_1=1$, then $\mx(W_1)=1$ and so
$U_{W}(v)=W_1$ for all $v\in W_1$ by Remark \ref{rc1} (i).
For $v_1\in W_1$, there exists $w_1\in W_2$ 
such that $e_1=v_1w_1\in E$.
Since $|J_{W_2}|>1$,
there exists $v_2\in W_2$ such that $w_1v_2\in E^{*}$,
and for this $v_2$, there exists $w_2\in W_1$ with $w_2\ne v_1$
such that $e_2=v_2w_2\in E$.
Certainly, $(e_1,e_2)$ is an R-cycle
because $w_2v_1\in E^{*}$.
In case of $m_1>1$ and $m_2>1$,
since $|W_i|\geq m_1+m_2>\mx(W_i)+1$ for $i=1,2$,
by virtue of Remark \ref{rc1} (iii),
$\UG_{W}$ satisfies the condition (UC).
Hence, by Lemma \ref{PS1},
$\Sy_{W}$ has an R-cycle
and so does $\Sy$.

Suppose next that $n>2$
and the statement holds for
all numbers between $2$ and $n-1$.
If $\UG_{W}$
satisfies the condition (UC),
then it has an R-cycle by Lemma \ref{PS1}.
We may assume therefore that $\UG_{W}$ fails to satisfy
the condition (UC); thus
there exists $i$ such that
$|U_{W_i}(v)\cap U_{W_i}(w)|\leq 1$ for some $v,w\in W_i$.
By Remark \ref{rc1}, for such $i$,
$\mx(W_i)=1$ if and only if $|W_i|=1$,
and it holds that either
\begin{equation}
\begin{array}{lll}
&|J_{W_i}|\leq 1\ \mbox{ and }\ |W_i|= \mx(W_i)\\
\mbox{or}
&|J_{W_i}|= 2\ \mbox{ and }\ |W_i|= \mx(W_i)+1.\\
\end{array}
\end{equation}

In case that there exists $i\in\{1,\cdots,n\}$
such that $|W_i|=\mx(W_i)$, say $i=n$,
we consider the R-subgraph $\Sy_{V'}$
with $V'=V\setminus V_n=V_1\cup\cdots\cup V_{n-1}$.
By Lemma \ref{SimSup} and the assumption of the statement, 
$$\begin{array}{lll}
&|I_{V'}(V')|-|I(V')|\leq |W_n|=\mx(W_n),\\
&|I(V')|=|I(V)|-|I(V_n)|\\
\mbox{and}& |I(V_n)|=|V_n|-|W_n|\geq (2m_n+1)-\mx(W_n),\\
\end{array}
$$
and so we have that
$$
\begin{array}{lll}
|I_{V'}(V')| &\leq |I(V')|+\mx(W_n) \\
&=|I(V)|-|I(V_n)|+\mx(W_n)\\
&\leq n-(2m_n+1)+2\mx(W_n)\leq n-1.
\end{array}
$$
By our induction hypothesis,
$\Sy_{V'}$ has an R-cycle.
We may assume therefore
that $|W_i|> \mx(W_i)$ for all $i\in\{1,\cdots,n\}$.
When this is the case, $m(W_i)>0$ for all $i$.

Then, by (3), there exists $i\in\{1,\cdots,n\}$
such that
$|W_i|= \mx(W_i)+1$ and $|J_{W_i}|= 2$; say $J_{W_i}=\{1,2\}$.
In addition, in this case, as mentioned at (3) above,
$\mx(W_i)> 1$ because of $|W_i|>1$.
We may here assume $i=n$.
Let $W_n=W_{n1}\cup W_{n2}$
with $W_{n1}=\{v_1\}$ and $W_{n2}=\{v_2,\cdots, v_{q+1}\}$,
where $q=\mx(W_n)$.
Since $|V_n|\geq 2m_n+1>|V_{n1}|+|V_{n2}|$,
there exists $k\in\{1,\cdots,l_{n}\}$
such that $k\ne 1,2$ and $I(V_n)\supseteq V_{nk}\ne\emptyset$,
where $l_n=|\cg(V_n)|$.
We may assume $k=3$.
Let $v_0\in V_{n3}$ and set $X_n=\{v_0,v_1,v_2\}$
and $V^{(1)}=V_1\cup\cdots\cup V_{n-1}\cup X_n$.
We should note $|V^{(1)}|<|V|$
because $|W_{n2}|=\mx(W_n)>1$.
We consider here the R-subgraph $\Sy_{V^{(1)}}$.
It is obvious that
$\UG_{X_n}$ satisfies (UC).
Let $X_{n1}=\{v_1\}$,
$X_{n2}=\{v_2\}$
and $X_{n3}=\{v_0\}$.
Then $\cg(X_n)=\{X_{n1},X_{n2},X_{n3}\}$
satisfies the R-colouring condition (RC)
because $X_{nj}\subseteq V_{nj}$
and $\{V_{n1},V_{n2},V_{n3}\}$ satisfies (RC).
Hence $\CG_n(V^{(1)})=
\{V_1,\cdots,V_{n-1},X_n\}$ is an R-colouring of $\Sy_{V^{(1)}}$.
We set $Y_n=V_n\setminus X_n$; thus $Y_n=V\setminus V^{(1)}$.
By Lemma \ref{SimSup},
$|I_{V^{(1)}}(V^{(1)})|\leq |I(V^{(1)})|+|Y_n|
-|I(Y_n)|.$
Since
$|I(V^{(1)})|=|I(V)|-(|I(V_n)|-1)$
and $|Y_n|-|I(Y_n)|
=|W_n|-2= \mx(W_n)-1$,
we have that
$$
\begin{array}{lll}
|I_{V^{(1)}}(V^{(1)})| & \leq |I(V)|-(|I(V_n)|-1)+(\mx(W_n)-1)\\
&\leq n-|I(V_n)|+\mx(W_n).\\
\end{array}$$
Moreover, because of $|V_n|\geq 2m_n+1$, we see that
$$|I(V_n)|=|V_n|-|W_n|\geq (2m_n+1)-(\mx(W_n)+1)\geq m_n,$$
and hence, 
$|I_{V^{(1)}}(V^{(1)})|\leq n-m_n+\mx(W_n)\leq n.$
In addition, $|X_n|\geq 2\mx(X_n)+1$,
in fact, $\mx(X_n)=1$ and $|X_n|=3= 2\mx(X_n)+1$.
That is,
$\Sy_{V^{(1)}}$ satisfy all of the conditions
supposed for $\Sy$ in the statement.
Let $W^{(1)}=W_1^{(1)}\cup\cdots\cup W_{n}^{(1)}$,
where $W_i^{(1)}=V_i\setminus I_{V^{(1)}}(V_i)$ $(i=1,\cdots,n-1)$
and $W_n^{(1)}=\{v_1,v_2\}$.
If 
$\UG_{W^{(1)}}$ fails to
the condition (UC)
and $|W_i^{(1)}|> \mx(W_i^{(1)})$ for all $i\in\{1,\cdots,n\}$,
then we can proceed with this procedure,
and get R-subgraphs $\Sy_{V^{(1)}},\Sy_{V^{(2)}},\cdots$.
On the other hand,
$|V|>|V^{(1)}|>|V^{(2)}|>\cdots$,
and therefore, there exists $p>0$
such that $\Sy_{V^{(p)}}$ satisfies either (UC) or
$|W_i^{(p)}|= \mx(W_i^{(p)})$
for some $i\in\{1,\cdots,n\}$.
In either case,
we have already seen that $\Sy_{V^{(p)}}$ has an R-cycle.
$\Box$
\vskip12pt

In the above theorem,
the assumption that
$|V_i|\geq 2\mx(V_i)+1$ for each $i\in\{1,\cdots, n\}$
cannot be dropped.
Let $\Sy=(V,E,E^{*})$ be the R-graph
which is described in Example \ref{ex1},
and let $\Sy'=(V',E',E^{*\prime})$ be the R-graph
with $V'=V\cup\{w_1,w_2,w_3,w_4\}$,
$E'=E$ and
$$E^{*\prime}=E^{*}\cup\{v_iw_1, v_jw_2, v_kw_3, v_lw_4\ |\
i=1,2,\ j=3,4,5,\ k=6,7,8,\ l=9,10 \}.$$
Then,
$V'_{11}=\{v_{1}\}$,
$V'_{12}=\{v_2\}$,
$V'_{13}=\{w_1\}$,
$V'_{21}=\{v_3,v_5\}$,
$V'_{22}=\{v_4\}$,
$V'_{23}=\{w_2\}$,
$V'_{31}=\{v_6,v_8\}$,
$V'_{32}=\{v_7\}$,
$V'_{33}=\{w_3\}$,
$V'_{41}=\{v_9\}$,
$V'_{42}=\{v_{10}\}$,
$V'_{43}=\{w_4\}$
and $V'_i=\bigcup_{j=1}^{3}V'_{ij}$.
That is,
$\Sy'$ is an R-colouring R-graph
with $\CG_4(V')=\{V'_1,V'_2,V'_3,V'_4\}$
and $\Gh^{*}_{V_i^{\prime}}=(V_i^{\prime},E^{*\prime}_{V_i^{\prime}})$
is a complete 3-partite graph.
In addition, $I(V')=\{w_1,w_2,w_3,w_4\}$ and so
$|I(V')|=4$.
Since $\mx(V_1)=\mx(V_4)=1$
and $\mx(V_2)=\mx(V_3)=2$,
we see that
$|V_i|=2\mx(V_i)+1$ for $i=1,4$
but $|V_i|=2\mx(V_i)$ for $i=2,3$.
As has been pointed out,
$\Sy$ is an R-colouring R-graph
with $\CG_4(V)$
and it has no R-cycles.
Hence $\Sy'$ has also no R-cycles,
because $V=V'\setminus I(V')$ in $\Sy'$ and
$\Sy'_{V}$ is isomorphic to $\Sy$.

Now,
let $\Sy$ be an R-colouring R-graph with
$\CG_n(V)=\{V_1,\cdots,V_n\}$.
By Remark \ref{rc1} (i),
if $\mx(V_i)=1$ for $i\in \{1,\cdots,n\}$,
$U(v)=V_i$ for all $v\in V_i$; thus
$\Gh^{*}_{V_i}=(V_i,E^{*}_{V_i})$
is a complete graph.
Therefore, if $\mx(V_i)=1$ for every $1\leq i\leq n$,
then $\Gh^{*}=(V,E^{*})$ is a disjoint union of complete graphs
and $\UG$ coincides with the set of components of $\Gh^{*}$;
thus $\UG=\CG_n(V)$.
In such case,
we define $\CG_n(V)$ to be a simple R-colouring of $\Sy$.
It is obvious that $\CG_n(V)$ is a simple R-colouring of $\Sy$
if and only if $\Sy$ is an R-graph satisfying the condition (SC):
\vskip6pt

\noindent
\mbox{{\rm (SC)}}\quad
$U,U'\in\UG\Longrightarrow \mbox{ either }U\cap U'=\emptyset
\mbox{ or } U=U'.$

\begin{df}\label{Rsimple}
Let $\Sy=(V,E,E^{*})$ be an R-graph.
$\Sy$ is R-simple
if the following conditions are satisfied$:$

\mbox{\rm (i)}
$\Sy$ has a simple R-colouring; thus
$\Sy$ satisfies \mbox{ {\rm (SC)}},

\mbox{\rm (ii)}
there exist no R-cycles of length $2$ consisting of edges;
if there exists $vw, v'w'\in E$
such that $vv'\in E^{*}$,
then $ww'\not\in E^{*}$.
\end{df}

Let $\Sy$ be an R-simple R-graph
with the set $\UG$ of R-neighbour sets.
By (i),
we can define the graph 
whose vertex set $\Vt=\UG$
and whose edge set $\Ed=\{UU'\ |\ U,U'\in \UG,$
there exist $v\in U$ and $v'\in U'$
such that $vv'\in E\}$; the graph $(\Vt, \Ed)$
is denoted by $\Sy/\UG$.
Then $|\Vt|=|\UG|$,
and for $U\in\Vt$, $d_{\Vts}(U)=\sum_{v\in U}d_{V}(v)$.
Moreover,
$\Sy/\UG$ has no loops by (R),
and has no multiple edges by (ii);
that is, $\Sy/\UG$ is a simple graph.
We call $\Sy/\UG$ the induced  simple graph of $\Sy$.

If $\Sy$ has an R-cycle then
it induces the cycle in $\Sy/\UG$.
Conversely,
if $\Sy/\UG$ has a cycle
then the origin of it in $\Sy$
is either an R-cycle
or a cycle in the base graph $\Gh=(V,E)$.
Hence we have

\begin{lemm}\label{RSP}
Let $\Sy$ be an R-simple R-graph with the base graph $\Gh=(V,E)$
and the set $\UG$ of R-neighbour sets.
Suppose that $\Gh$ has no cycles.
Then $\Sy$ has an R-cycle if and only if $\Sy/\UG$ has a cycle.
\end{lemm}

\begin{df}\label{Rcon}
Let $\Sy$ be an R-simple R-graph
with the set $\UG$ of R-neighbour sets.
Then $U$ and $U'$ in $\UG$ are said to be R-connected
if there exists a finite sequence $U_0e_1U_1\cdots e_pU_p$
whose terms are alternately R-neighbour sets $U_q$'s
and edges $e_q$'s in $E$ with $e_q=v_qw_q$ such that
$U_0=U$, $U_p=U'$, $v_q\in U_{q-1}$, $w_q\in U_q$.
\end{df}

In an R-simple R-graph,
'R-connected' means simply 'connected' in the induced simple graph of it.
Since R-connection is an equivalence relation on $\UG$,
there exists a decomposition of $\UG$ into subsets $\UG_i$'s
$(1\leq i\leq m)$ for some $m>0$ such that $U,U'\in \UG$
are R-connected if and only if both $U$ and $U'$ belong to
the same set $\UG_i$.
The subgraphs $\Sy_{W_1},\cdots,\Sy_{W_m}$ of $\Sy$ generated by $W_i$'s
are called the R-components of $\Sy$,
provided $W_i=\bigcup_{U\in\UG_i}U$ for each $i\in\{1,\cdots,m\}$
and $V=\bigcup_{i=1}^{m}W_i$.
If $\Sy$ has exactly one R-component
then $\Sy$ is R-connected.

Recall that $\NGG(\Sy)=\{U\in \UG\ |\ |U|=1\}$.
We set
$\LG(\Sy)=\{U\in \UG\ |\ |U|>2\}$
and $\MG(\Sy)=\{U\in \UG\ |\ |U|=2\}$.

\begin{lemm}\label{ConnSim}
Let $\Sy=(V,E,E^{*})$ be an R-simple R-graph
with the set $\UG$ of R-neighbour sets.
Suppose that $\Sy$ is R-connected.
Then $\Sy$ has an R-cycle
if and only if $|V|-|\UG|-\omega+1\ne 0$,
where $\omega$ is the number of components of $\Gh$.

In particular,
if $\Sy$ is proper and $|\LG(\Sy)|\geq |\NGG(\Sy)|$
then $\Sy$ has an R-cycle.
\end{lemm}

\noindent
{\bf Proof.}\quad
Let $\Gh_i=(V_i, E_i)$ $(i=1,\cdots, \omega)$
be the components of the base graph $\Gh=(V,E)$.
Since $\Gh_i$ is connected,
there exists a spanning tree $\Gh_i^{\prime}=(V_i,E_i^{\prime})$ of $\Gh_i$.
Let $E^{\prime}=\bigcup_{i=1}^{\omega}E_i^{\prime}$,
$\Gh^{\prime}=(V,E^{\prime})$
and $\Sy^{\prime}=(V,E^{\prime},E^{*})$.
It is obvious that $\Sy^{\prime}$ is R-simple and R-connected.
Moreover, $\Sy$ has an R-cycle if and only if so does $\Sy^{\prime}$.
Since $\Gh^{\prime}$ has no cycles,
by virtue of Lemma \ref{RSP},
$\Sy^{\prime}$ has an R-cycle if and only if the induced simple graph
$\Sy^{\prime}/\UG=(\Vt^{\prime}, \Ed^{\prime})$ has a cycle.

On the other hand,
since $\Sy^{\prime}$ is R-connected,
$\Sy^{\prime}/\UG$ is connected,
and so $\Sy^{\prime}/\UG$ is a tree if and only if
\begin{equation}
|\Ed^{\prime}|=|\Vt^{\prime}|-1.
\end{equation}
We set $D=\sum_{U\in\Vts^{\prime}}d_{\Vts^{\prime}}(U)$
and $d_i=\sum_{v\in V_i}d_{V}(v)$,
where $d_{V}(v)$ means the degree of $v$ not in $\Gh$ but in $\Gh^{\prime}$.
Recall that $d_i=2|V_i|-2$ because $\Gh_i^{\prime}$ is a tree.
Hence we have
$$D
=\sum_{i=1}^{\omega}d_i
=\sum_{i=1}^{\omega}(2|V_i|-2)=2|V|-2\omega.$$
Since $|\Vt^{\prime}|=|\UG|$ and
generally $D=2|\Ed^{\prime}|$, the condition (4) can be replaced
by $|V|-|\UG|-\omega+1= 0$.
That is, $\Sy^{\prime}/\UG$ has a cycle if and only if
$|V|-|\UG|-\omega+1\ne 0$.

Now,
if $\Sy$ is proper,
$2|V|-2\omega=D\geq 3|\LG(\Sy)|+2|\MG(\Sy)|+|\NGG(\Sy)|$
and $|\UG|=|\LG(\Sy)|+|\MG(\Sy)|+|\NGG(\Sy)|$.
Hence, in particular,
if $|\LG(\Sy)|\geq |\NGG(\Sy)|$, then
$$
|V|-|\UG|-\omega+1 \geq \frac{1}{2}(|\LG(\Sy)|-|\NGG(\Sy)|)+1>0,
$$
and so $\Sy$ has an R-cycle.
$\Box$
\vskip12pt

\begin{theom}\label{PS2}
Let $\Sy=(V,E,E^{*})$ be an R-simple R-graph
with the set $\UG$ of R-neighbour sets.
Then
$\Sy$ has an R-cycle
if and only if there exists an R-component
$\Sy_{W}=(W,E_{W},E^{*}_{W})$ of $\Sy$
with the set $\UG_{W}$ of R-neighbour sets
such that $|W|-|\UG_{W}|-\omega+1\ne 0$,
where $\omega$ is the number of
components of $\Gh_{W}=(W,E_{W})$.

In particular,
if $\Sy$ is proper and $|\LG(\Sy)|\geq |\NGG(\Sy)|$
then $\Sy$ has an R-cycle.
\end{theom}

\noindent
{\bf Proof.}\quad
If there exists an R-cycle in $\Sy$,
then it exists in an R-component.
Hence $\Sy$ has an R-cycle
if and only if there exists an R-component
$\Sy_{W}$
of $\Sy$
which has an R-cycle.
Let $\Sy_{W}=(W,E_{W},E^{*}_{W})$
be an R-component of $\Sy$ with the set $\UG_{W}$ of R-neighbour sets,
and let $\omega$ be the number of components
of $\Gh_{W}=(W,E_{W})$.
By Lemma \ref{ConnSim},
$\Sy_{W}$ has an R-cycle if and only if
$|W|-|\UG_{W}|-\omega+1\ne 0$.

Now,
if $|\LG(\Sy)|\geq |\NGG(\Sy)|$
then there exists an R-component $\Sy_{W}$
such that $|\LG(\Sy_{W})|\geq |\NGG(\Sy_{W})|$,
and in addition, if $\Sy$ is proper then so is $\Sy_{W}$.
Hence, in this case, $\Sy_{W}$ has an R-cycle by Lemma \ref{ConnSim}.
This completes the proof.
$\Box$
\vskip12pt

\section{PROOF OF THEOREMS}
\label{profth}

In what follows,
let $G$ be a non-abelian locally free group in which
there exists a free subgroup $H$ with basis $X$
such that $|H|=|G|$.
Since  $G$ is non-abelian,
we may assume that $H$ is non-abelian.
If the rank of $H$ is infinite then $|X|=|H|$.
On the other hand,
if the rank of $H$ is finite then
there exists a free subgroup of $H$
whose rank is countable; for instance,
the derived subgroup $[H,H]$ of $H$
is a free group of countable rank.
Therefore, we may assume here that $|X|=|G|$. 
Let $R$ $(\ni 1)$ be a ring with no zero divisors.
We suppose that $|R|\leq |G|$. 
Since $|X|=|G|\geq \aleph_{0}$,
we can divide $X$ into three subsets $X_1$, $X_2$ and $X_3$
each of whose cardinality is $|X|$.
Let $\sigma_i$ be a bijection from $X$ to $X_i$ $(i=1,2,3)$.
For $ x\in X$,
$x^{(i)}$ denote the image of $x$ by $\sigma_i$.
Let 
$RG$ be the group ring of $G$ over $R$.
Since $|RG|=|X|$,
there exists a bijection $\psi$ from $X$ to $RG\setminus \{0\}$.
Let $x\in X$, and let
\begin{equation}
\psi(x)=\sum_{i=1}^{m_x}\alpha_{xi}f_{xi},\
\mbox{ where }\ \alpha_{xi}\in R,\
f_{xi}\in G,\
m_x>0,
\end{equation}
each of them depends on $x$, and
%
$f_{xi}\ne f_{xj} \mbox{ if } i\ne j.$
%
Since $G$ is locally free,
for the subset $\{x^{(1)},x^{(2)},x^{(3)}, 
f_{xi}\ |\ 1\leq i\leq m_x \}$ of $G$,
there exist elements
$z_{x1}, z_{x2}, z_{x3}$
in $G$ which satisfy the assertions of Lemma \ref{fg_b} (2).
We define $\varepsilon(x)$ by
\begin{equation}
\begin{array}{lll}
\varepsilon(x)
&\displaystyle
=\sum_{k=1}^{3}\sum_{l=1}^{3}x^{(k)}z_{xl}\psi(x)z_{xl}+1\\
&\displaystyle
=\sum_{k=1}^{3}\sum_{l=1}^{3}\sum_{i=1}^{m_x}
\alpha_{xi}x^{(k)}z_{xl}f_{xi}z_{xl}+1.\\
\end{array}
\end{equation}
Let $\xi_{x}(k,l,i)
=x^{(k)}z_{xl}f_{xi}z_{xl}$.
Then the assertions of Lemma \ref{fg_b} (2) mean the the following:

\begin{rema}\label{Rem1}

$\mbox{{\rm (i)}}$
$\xi_{x}(k,l,i)=\xi_{x}(h,n,j)$
if and only if $(k,l,i)=(h,n,j)$.

$\mbox{{\rm (ii)}}$
Let $p>0$, $1\leq l_{q},n_{q}\leq 3$ and
$1\leq i_{q},j_{q}\leq m_x$,
where $1\leq q\leq p$.
If $\prod_{q=1}^{p}\xi_{x}(1,l_{q},i_{q})^{-1}\xi_{x}(1,n_{q},j_{q})=1$,
then 
either $n_{q}=l_{q+1}$ for some $q\in\{1,\cdots,p-1\}$
or $(l_{q},i_{q})=(n_{q},j_{q})$
for some $q\in\{1,\cdots,p\}$.

\end{rema}

Let $\rho$ be the right ideal of $RG$
generated by all $\varepsilon(x)$'s,
that is,
\begin{equation}
\rho=\sum_{x\in X}\varepsilon(x)RG.
\end{equation}
Let $r=\sum_{t=1}^{m}r_{t}$ be a non-zero element of $\rho$,
where $0\ne r_{t}\in\varepsilon(x_{t})RG$
with $x_{t}\in X$.
Then there exist $n_{t}>0$,
$\beta_{tj}\in R$ with $\beta_{tj}\ne 0$
and $g_{tj}\in G$
such that
$r_{t}=\varepsilon(x_{t})\sum_{j=1}^{n_{t}}\beta_{tj}g_{tj}$
with $g_{tj}\ne g_{ti}$ $(j\ne i)$.
In what follows,
we simply write $m_{t}$, $\alpha_{ti}$, $x_{tk}$, $z_{tl}$ and $f_{ti}$ instead of
$m_{x_t}$, $\alpha_{x_ti}$, $x_t^{(k)}$, $z_{x_tl}$ and $f_{x_ti}$, respectively.
By the expression of $\varepsilon(x_{t})$ as in (6),
we have that
\begin{equation}
r_{t}=\sum_{k,l=1}^{3}\sum_{i=1}^{m_t}\sum_{j=1}^{n_t}
\alpha_{ti}\beta_{tj}x_{tk}z_{tl}f_{ti}z_{tl}
g_{tj}+\sum_{j=1}^{n_t}\beta_{tj}g_{tj}, \mbox{ where } m_t,n_t>0.
\end{equation}
To prove Theorem \ref{MainTh},
by virtue of Proposition \ref{ForProp},
we shall show that
$\rho$ is a proper right ideal of $RG$.
By making use of
graph-theoretic results obtained in the previous section,
we shall prove $r=\sum_{t=1}^{m}r_{t}\ne 1$.
To connect the problem with our graphical method,
we prepare the following notations.

For $1\leq t\leq m$,
let $n_{t}$ and $m_{t}$ be as described in (8).
we set
\begin{equation}
\begin{array}{lll}
&P_t=\{(t,k,l,i,j)\ |\ 1\leq k,l\leq 3,\
1\leq i\leq m_{t},\ 1\leq j\leq n_{t}\}\\
\mbox{ and }\
&Q_t=\{(t,j)\ |\ 1\leq j\leq n_{t}\}.\\
\end{array}
\end{equation}
Moreover,
for $v=(t,k,l,i,j)\in P_t$ and $w=(t,j)\in Q_t$,
we set
\begin{equation}
\eta(v)
=x_{tk}z_{tl}f_{ti}z_{tl}g_{tj}
\ \mbox{ and }\
\eta(w)=g_{tj}.
\end{equation}
Then we can replace the expression (8) of $r_{t}$
by the following expression:
\begin{equation}
r_{t}=\sum_{v\in P_{t}}\gamma_{v}\eta(v)
+\sum_{w\in Q_{t}}\gamma_{w}\eta(w),
\mbox{ where }
\gamma_{v}=\alpha_{ti}\beta_{tj}
\mbox{ and } \gamma_{w}=\beta_{tj}.
\end{equation}

Let $P=\bigcup_{t=1}^{m}P_{t}$ and $Q=\bigcup_{t=1}^{m}Q_{t}$.
We regard $W=P\cup Q$
as the set of vertices
and $E=\{vw\ |\ v,w\in W, v\ne w \mbox{ and } \eta(v)=\eta(w)\}$
as the set of edges,
and consider the R-graph $\Sy=(W,E,E^{*})$,
where $vv'\in E^{*}$ if and only if $v, v'\in P_{t}$
such that $v=(t,k,l,i,j)$ and $v'=(t,k',l,i,j)$ with $k\ne k'$; thus
$U(v)=\{(t,k',l,i,j)\ |\ 1\leq k'\leq 3\}$ for $v=(t,k,l,i,j)\in P_t$,
$U(w)=\{w\}$ for $w\in Q_t$ and $\UG=\{U(v)\ |\ v\in W\}$
(in the proof of Theorem \ref{MainTh},
in fact,
the above vertices set $W=P\cup Q$ is replaced by
$V=P^{*}\cup Q^{*}$; see below for detail).
We shall then show that
there exist some isolated vertices in $\Sy$
which make $r=1$ false.
To do this,
by making use of Theorem \ref{PS3},
we shall first show that
there exists a suitable number of isolated vertices
in the subgraph of $\Gh=(W,E)$ generated by $P_{t1}$
in Lemma \ref{Mt} after preparing two remarks,
where $P_{t1}=\{v\ |\ v=(t,1,l,i,j)\in P_t\}$.

Let
$M_t=\{(l,i,j)\ |\  1\leq l\leq 3,\
1\leq i\leq m_{t},\ 1\leq j\leq n_{t}\}.$
For $v=(t,k,l,i,j)\in P_{t}$ and $\mu=(l,i,j)\in M_{t}$,
we write $v=(t,k,\mu)$.
Let $\mu=(l,i,j)$ and $\mu'=(l',i',j')$ be in $M_{t}$.
If $j=j'$ then
$\eta(t,k,\mu)=\eta(t,k',\mu')$ if and only if
$(k,l,i)=(k',l',i')$ by Remark \ref{Rem1} (i).
If $(l,i)=(l',i')$,
then $\eta(t,1,\mu)=\eta(t,1,\mu')$ implies
$g_{tj}=g_{tj'}$,
and so $j=j'$; thus $\mu=\mu'$.
Therefore,
if $\mu\ne\mu'$,
then $\eta(t,1,\mu)=\eta(t,1,\mu')$
implies
$j\ne j'$ and $(l,i)\ne (l',i')$.
Moreover, it is obvious that
$\eta(t,1,\mu)=\eta(t,1,\mu')$ holds
if and only if
$\eta(t,k,\mu)=\eta(t,k,\mu')$ holds for all $k\in\{1,2,3\}$.
Hence we have

\begin{rema}\label{Rem3}
Let $t\in\{1,\cdots,m\}$,
and let $v,v'\in P_t$ with $v=(t,k,\mu)$ and $v'=(t,k',\mu')$,
where $\mu=(l,i,j)$ and $\mu'=(l',i',j')$.

\mbox{{\rm (i)}}\
Suppose that $\eta(v)=\eta(v')$.
If $j=j'$,
then $v=v'$.

\mbox{{\rm (ii)}}\
Suppose that $\eta(v)=\eta(v')$.
If $k=k'$ and either $j=j'$ or $(l,i)=(l',i')$,
then $\mu=\mu'$.

\mbox{{\rm (iii)}}\
$\eta(t,k,\mu)=\eta(t,k,\mu')$ holds for some $k\in\{1,2,3\}$
if and only if it holds for any $k\in\{1,2,3\}$.
\end{rema}

For $\mu$ and $\mu'$ in $M_{t}$,
define the relation $\mu\sim \mu'$
by $\eta(t,1,\mu)=\eta(t,1,\mu')$.
It is obvious that $\sim$ is a equivalence relation on $M_{t}$.
Let $C_{t}(\mu)$ be the equivalence class of $\mu\in M_{t}$
and let $N_{t}=\{\mu\in M_{t}\ |\ C_{t}(\mu)=\{\mu\}\}$.
Since $N_{t}\subseteq M_{t}$ and $\nu\not\sim\nu'$
for $\nu,\nu'\in N_{t}$ with $\nu\ne\nu'$,
as a complete set of representatives for $M_{t}/\sim$,
we can choose a set $T_{t}$ which satisfies $N_{t}\subseteq T_{t}$.
For $\mu=(l,i,j)\in T_{t}$,
let $\gamma_{(t,1,\mu)}=\alpha_{ti}\beta_{tj}$,
where $\alpha_{ti}\beta_{tj}$ is as described in (11).
We set
$$
\gamma_{(t,1,\mu)}^{*}=\sum_{\mu'\in C_{t}(\mu)}\gamma_{(t,1,\mu')}\
\mbox{ and }\ M_{t}^{*}=\{\mu\in T_{t}\ |\ \gamma_{(t,1,\mu)}^{*}\ne 0\}.
$$
By the definition of $M_{t}^{*}$ and Remark \ref{Rem3} (i) (iii),
we have

\begin{rema}\label{Rem4}
Let $t\in \{1,\cdots m\}$ and $k,k'\in\{1,2,3\}$.
For $\mu,\mu'\in M_{t}^{*}$,
suppose $\eta(t,k,\mu)=\eta(t,k',\mu')$.
Then $k=k'$ if and only if $\mu=\mu'$.
\end{rema}

Now in (11),
replacing $P_{t}$ by
$P_{t}^{*}=\{(t,k,\mu)\ |\ \mu\in M_{t}^{*},\ 1\leq k\leq 3\}$,
we can use the following expression of $r_{t}$:
\begin{equation}
r_{t}=\sum_{v\in P_{t}^{*}}\gamma_{v}^{*}\eta(v)
+\sum_{w\in Q_{t}}\gamma_{w}\eta(w).
\end{equation}
%


\begin{lemm}\label{Mt}
Let $n_{t}=|Q_t|$ be as described in $(9)$,
and let $M_{t}^{*}$ and $N_{t}$ as above.
Then $|N_{t}|>n_{t}$ for all $t\in\{1,\cdots,m\}$.
In particular, $|M_{t}^{*}|> n_{t}$.
\end{lemm}

\noindent
{\bf Proof.}\quad
Since $M_{t}^{*}\supseteq N_{t}$,
it suffices to show that $|N_{t}|> n_{t}$.
Suppose, to the contrary,
that $|N_{t}|\leq n_{t}$ for some $t\in\{1,\cdots,m\}$.
If $n_{t}=1$ then $|N_{t}|=|M_{t}|$ by Remark \ref{Rem3} (i),
which implies $|N_{t}|=3m_{t}\geq 3>1=n_{t}$, a contradiction.
Hence we have $n_{t}>1$.

Let $\Gh=(V,E)$ be the graph with the vertex set $V=M_t$
and the edge set $E$ defined by
$$\{vw\ |\ v,w\in V,\ v\ne w,\ \eta(t,1,v)=\eta(t,1,w)\}.$$
It is obvious that $\Gh_{W}=(W,E)$ is a clique graph,
where $W=V\setminus N_{t}$.
We set $V_{jl}=\{(l,i,j)\in V\ |\ 1\leq i\leq m_{t}\}$
and $V_j=\bigcup_{l=1}^{3}V_{jl}$; thus $V=\bigcup_{j=1}^{n_t}V_j$.
Let $\Gh^{*}_{V_j}=(V_j,E^{*}_{V_j})$ $(j=1,\cdots, n_t)$
be the complete 3-partite graph
with the partite set $\{V_{j1},V_{j2},V_{j3}\}$
and let $\Gh^{*}=(V,E^{*})=\bigcup_{j=1}^{n_t}\Gh^{*}_{V_j}$.
Since $U(v)=V_j\setminus V_{jl}\cup\{v\}$ for $v\in V_{jl}$
and $\Gh_{W}$ is a clique graph,
by Remark \ref{Rem3} (i),
$\Sy=(V,E,E^{*})$ satisfies (R), and so $\Sy$
is an R-graph; in fact, it is a non-empty clique R-graph.
In addition, since $\Gh^{*}_{V_j}$ is the complete 3-partite graph,
$\Sy$ is an R-colouring R-graph
with $\CG_{n_t}(V)=\{V_1,\cdots,V_{n_t}\}$.
Since $m_t>0$, we see then that $|V_j|=3m_t\geq 2m_t+1=2\mx(V_j)+1$
for each $j\in\{1,\cdots,n_t\}$.
Moreover,
according to our hypothesis,
$|N_{t}|\leq n_{t}$, that is,
$|I(V)|\leq n_t$.
Hence, by virtue of Theorem \ref{PS3}, a clique R-graph
$\Sy$ has an R-cycle consisting of edges.
That is,
there exist $p>1$ and edges 
$e_1,\cdots,e_p\in E$ with $e_{q}=v_{q}w_{q}$ $(1\leq q\leq p)$
such that
all of $v_q$'s and $w_q$'s are different from each other,
$w_{q}v_{q+1}\in E^{*}$ $(1\leq q\leq p-1)$
and $w_{p}v_{1}\in E^{*}$.
Let $v_{q}=(l_{q},i_{q},j_{q})$
and
$w_{q}=(l_{q}^{\prime},i_{q}^{\prime},j_{q}^{\prime})$,
where $1\leq q\leq p$.
Let $\xi_{t}(1,l_q,i_q)
=x_{t1}z_{tl_q}f_{ti_q}z_{tl_q}$.
Then $e_q=v_qw_q\in E$ implies
$$\xi_{t}(1,l_{q},i_{q})g_{tj_{q}}
=\eta(t,1,v_{q})
=\eta(t,1,w_{q})
=\xi_{t}(1,l_{q}^{\prime},i_{q}^{\prime})g_{tj_{q}^{\prime}}.$$
Moreover,
$w_qv_{q+1}\in E^{*}$ and $w_pv_1\in E^{*}$ mean that
$j_{q}^{\prime}=j_{q+1}$, $j_{p}^{\prime}=j_{1}$,
and $l_{q}^{\prime}\ne l_{q+1}$.
Hence we have
$$\prod_{q=1}^{p}\xi_{t}(1,l_{q},i_{q})^{-1}
\xi_{t}(1,l_{q}^{\prime},i_{q}^{\prime})=1\
\mbox{ with }\ l_{q}^{\prime}\ne l_{q+1}
\ (1\leq q\leq p-1).$$
Since $v_q\ne w_q$ and $\eta(t,1,v_{q})=\eta(t,1,w_{q})$
for all $1\leq q\leq p$,
it follows from Remark \ref{Rem3} (ii)
that $(l_{q},i_{q})\ne (l_{q}^{\prime},i_{q}^{\prime})$ for all $1\leq q\leq p$.
However, this contradicts the assertion of Remark \ref{Rem1} (ii).
$\Box$
\vskip12pt

We are now in a position to prove Theorem \ref{MainTh}.

\noindent
{\bf Proof.}\quad
[Proof of Theorem \ref{MainTh}]
(1):
Let $\rho$ be as described in $(7)$;
a non-trivial right ideal of $RG$.
By virtue of Proposition \ref{ForProp},
it suffices to show that
$\rho$ is proper.
Let $r=\sum_{t=1}^{m}r_{t}$ be as described in $(11)$;
a non-zero element of $\rho$.
For $t\in\{1,\cdots,m\}$,
recall
$P_{t}^{*}=\{(t,k,\mu)\ |\ \mu\in M_{t}^{*},\ 1\leq k\leq 3\}$
and $Q_t=\{(t,j)\ |\ 1\leq j\leq n_{t}\}.$
We set $P^{*}=\cup_{t=1}^{m}P_{t}^{*}$
and $Q^{*}=(\cup_{t=1}^{m}Q_{t})\cup \{w_0\}$,
where $w_{0}=(0,0)$. We define $\gamma_{w_{0}}=-1$ $(\in R)$
and $\eta(w_{0})=1$
$(\in G)$ respectively.
Then $r=\sum_{v\in P^{*}}\gamma_{v}^{*}\eta(v)
+\sum_{w\in Q^{*}}\gamma_{w}\eta(w)+1$ by (12).
In order to prove that $\rho$ is proper,
it suffices to show that $r\ne 1$.
Suppose, to the contrary,
that $r-1=0$, that is,

\begin{equation}
\sum_{v\in P^{*}}\gamma_{v}^{*}\eta(v)
+\sum_{w\in Q^{*}}\gamma_{w}\eta(w)=0.
\end{equation}

Now,
we set $V=P^{*}\cup Q^{*}$
and let $\Gh=(V,E)$ be the graph
whose vertices are the elements of $V$
and whose edge set $E$
is defined as
$$E=\{vw\ |\ v,w\in V,\ v\ne w,\ \eta(v)=\eta(w) \mbox{ in } G\}.$$
By (13), $\Gh$ is proper and it is a non-empty clique graph.
Let $E^{*}=\{ vw\ |\ v=(t,k,\mu), w=(t,k',\mu)\in P^{*}, k\ne k'\}$.
By Remark \ref{Rem4}, $\Sy=(V,E,E^{*})$ satisfies (R$'$), and so
$\Sy$ is an R-graph,
and in fact, it is a non-empty clique R-graph.
We shall show that $\Sy$ has an R-cycle.
If $v\in V$,
$$
U(v)=\left\{\begin{array}{ll}
\{(t,k',l,i,j)\ |\ 1\leq k'\leq 3\} \mbox{ if } v=(t,k,l,i,j)\in P^{*}\\
\{v\} \mbox{ if } v\in Q^{*}\\
\end{array}\right. ,$$
and so $\UG=\{U(v)\ |\ v\in V\}$ satisfies the condition (SC).
Hence, either $\Sy$ is R-simple
or it has an R-cycle of length $2$ consisting of edges.
We may assume, therefore, that $\Sy$ is R-simple.
If $v\in P^{*}$ then $|U(v)|=3$,
and if $v\in Q^{*}$ then $|U(v)|=1$.
This means that
$\LG(\Sy)=\{U(v)\ |\ v\in P^{*}\}$
and $\NGG(\Sy)=\{U(v)\ |\ v\in Q^{*}\}$.
By Lemma \ref{Mt},
$|M_{t}^{*}|> n_{t}$, which implies
$$\textstyle |\LG(\Sy)|=(1/3)|P^{*}|
=\sum_{t=1}^{m}|M_{t}^{*}|>\sum_{t=1}^{m}n_{t}.$$
On the other hand,
$$\textstyle |\NGG(\Sy)|=|Q^{*}|=
\sum_{t=1}^{m}|Q_{t}|+1=\sum_{t=1}^{m}n_{t}+1,$$
and so $|\LG(\Sy)|\geq |\NGG(\Sy)|$.
Hence, by Theorem \ref{PS2},
a clique R-graph $\Sy$ 
has an R-cycle consisting of edges,
as desired.
By the definition of an R-cycle,
there exist $p>1$ and
$v_{q}=(t_{q},k_{q},\mu_{q}),
w_{q}=(t_{q},h_{q},\mu_{q})
\in P^{*}$
with $\mu_{q}=(l_{q},i_{q},j_{q})\in M_{t_{q}}^{*}$
$(1\leq q\leq p)$
such that all of $v_q$'s and $w_q$'s are different from each other,
$w_qv_q\in E^{*}$,
$e_q=v_qw_{q+1}\in E$ $(1\leq q\leq p-1)$
and $v_pw_1\in E$.
Hence, we have
\begin{equation}
\quad\quad\ \eta(v_{q})=\eta(w_{q+1})\ (1\leq q\leq p-1),\quad \eta(v_{p})=\eta(w_{1})
\end{equation}
\begin{equation}
\mbox{ and }\ k_{q}\ne h_{q}\ (1\leq q\leq p).
\end{equation}
For the sake of simplicity of notation,
the subscript $p+1$ means the subscript $1$;
we set $w_{p+1}=w_1$, $t_{p+1}=t_1$, $j_{p+1}=j_1$ $\cdots$.
Since $v_q\ne w_{q+1}$,
by Remark \ref{Rem3} (i),
$\eta(v_{q})=\eta(w_{q+1})$ implies
\begin{equation}
\mbox{ either }\ t_{q}\ne t_{q+1}\
\mbox{ or }\ j_{q}\ne j_{q+1}.
\end{equation}
By (10), the definition of $\eta(v_q)$,
$$\eta(v_{q})=x_{t_{q}k_{q}}\zeta_q\
\mbox{ and }\
\eta(w_{q+1})=x_{t_{q+1}h_{q+1}}\zeta_{q+1},
$$
where
$\zeta_q=z_{t_{q}l_{q}}f_{t_{q}i_{q}}
z_{t_{q}l_{q}}g_{t_{q}j_{q}}$,
and so
(14) implies that
\begin{equation}
\prod_{q=1}^{p}(x_{t_{q}k_{q}})^{-1}x_{t_{q+1}h_{q+1}}=1.
\end{equation}
Recall that $x_{tk}$'s are elements in $X$ which is a basis of a free group
and that $x_{tk}=x_{t'k'}$ if and only if $(t,k)=(t',k')$.
Now, in (17), if $t_{q}\ne t_{q+1}$,
it is obvious that
$x_{t_{q}k_{q}}\ne x_{t_{q+1}h_{q+1}}$.
If $t_{q}=t_{q+1}$ then $j_{q}\ne j_{q+1}$ because of (16),
and so $\mu_{q}\ne\mu_{q+1}$.
Since $\eta(v_{q})=\eta(w_{q+1})$ by (14),
Remark \ref{Rem4} implies $k_{q}\ne h_{q+1}$,
and hence,
$x_{t_{q}k_{q}}\ne x_{t_{q+1}h_{q+1}}$ again.
Moreover,
we have that
$x_{t_{q+1}h_{q+1}}\ne x_{t_{q+1}k_{q+1}}$
by (15).
Therefore it follows a contradiction
that $\prod_{q=1}^{p}(x_{t_{q}k_{q}})^{-1}x_{t_{q+1}h_{q+1}}\ne 1$. 
This complete the proof of (1).

(2):
If $K'$ is the prime field of $K$
then $|K'|\leq |G|$,
and therefore $K'G$ is primitive by (1).
Since $\Delta(G)=1$ by Lemma \ref{fg_b} (1),
the conclusion follows from Lemma \ref{Ps}.
$\Box$
\vskip12pt

Now, a ring $R$ is called a (right) strongly prime ring
if for each $0\ne \alpha\in R$,
there exists a finite subset $S(\alpha)$ of $R$
such that $\alpha S(\alpha)\beta\ne 0$
for all non-zero $\beta\in R$.
$S(\alpha)$ is called a (right) insulator of $\alpha$.
For instance, domains and simple rings are strongly prime.
Formanek's result \cite[Theorem]{For} on primitivity of
$RG$ for a domain $R$ was generalized to one for a strongly prime ring $R$
by Lawrence \cite{Law}.
The same situation holds for the case of our theorem.

\begin{coro}\label{subMainCor}
The assertion of Theorem \ref{MainTh} (1)
holds also for a strongly prime ring $R$.
\end{coro}

\noindent
{\bf Proof.}\quad
Let $\varphi(x)=\sum_{i=1}^{m_x}\alpha_{xi}f_{xi}$ $(x\in X)$
be as described in (5)
and $S(\alpha_{x1})$
$=\{\delta_{xq}\ |\ 1\leq q\leq d_x\}$ 
a right insulator of $\alpha_{x1}$.
Going back to the beginning of this section,
for $\{x^{(1)},x^{(2)},x^{(3)}, f_{xi}\ |\ 1\leq i\leq m_x\}$,
there exist elements $z_{xlq}$ $(1\leq l\leq 3, 1\leq q\leq d_x)$
which satisfy assertion of Lemma \ref{fg_b} (2).
We replace (6) by
$$
\begin{array}{lll}
\varepsilon(x)
&\displaystyle
=\sum_{k=1}^{3}\sum_{l=1}^{3}\sum_{q=1}^{d_x}x^{(k)}z_{xlq}\psi(x)\delta_{xq}z_{xlq}+1\\
&\displaystyle
=\sum_{k=1}^{3}\sum_{l=1}^{3}\sum_{q=1}^{d_x}\sum_{i=1}^{m_x}
\alpha_{xi}\delta_{xq}x^{(k)}z_{xlq}f_{xi}z_{xlq}+1.\\
\end{array}
$$
Then (8) is replaced by
\begin{equation}
\begin{array}{lll}
r_{t}=
&\displaystyle\sum_{k,l=1}^{3}\sum_{q=1}^{d_t}\sum_{i=1}^{m_t}\sum_{j=1}^{n_t}
\alpha_{ti}\delta_{tq}\beta_{tj}x_{tk}z_{tlq}f_{ti}z_{tlq}g_{tj}\\
&\hskip3cm \displaystyle+\sum_{j=1}^{n_t}\beta_{tj}g_{tj}, \mbox{ where } m_t,n_t,d_t>0.\\
\end{array}
\end{equation}
Let $A_{tj}=\{q\ |\ 1\leq q\leq d_x, 
\alpha_{ti}\delta_{tq}\beta_{tj}\ne 0 \mbox{ for some } i\}$.
Since $S(\alpha_{t1})$
is a right insulator of $\alpha_{t1}$,
for each $j\in\{1,\cdots,n_t\}$,
there exists $q\in\{1,\cdots,d_t\}$
such that $\alpha_{t1}\delta_{tq}\beta_{tj}\ne 0$,
and so $A_{tj}\ne\emptyset$.
For $q\in A_{tj}$,
let $B_{tj}(q)=\{i\ |\ 1\leq i\leq m_t,
\alpha_{ti}\delta_{tq}\beta_{tj}\ne 0\}$.
It is obvious that $B_{tj}(q)\ne\emptyset$.
The non-zero parts of (18) is here
replaced by the following expression:
\begin{equation}
\begin{array}{lll}
r_{t}=
&\displaystyle\sum_{k,l=1}^{3}
\sum_{j=1}^{n_t}\sum_{q\in A_{tj}}\sum_{i\in B_{tj}(q)}
\alpha_{ti}\delta_{tq}\beta_{tj}x_{tk}z_{tlq}f_{ti}z_{tlq}g_{tj}\\
&\displaystyle+\sum_{j=1}^{n_t}\beta_{tj}g_{tj},
\mbox{ where } n_t>0, |A_{tj}|>0, |B_{tj}(q)|>0.\\
\end{array}
\end{equation}
After this,
we renumber the elements in $A_{tj}$ and $B_{tj}(q)$,
and we can then follow the same proof as in Theorem \ref{MainTh} (1).

We can summarize the procedure as follows:
Let $A_{tj}=\{q_1,\cdots, q_{a_{tj}}\}$; $a_{tj}=|A_{tj}|$,
and for $q_s\in A_{tj}$, let $B_{tj}(q_s)=\{p_1,\cdots,p_{m_{tjs}}\}$;
$m_{tjs}=|B_{tj}(q_s)|$.
We set $A^{*}_{tj}=\{1,\cdots,3a_{tj}\}$
and $B^{*}_{tjl}=\{1,\cdots,m_{tjs}\}$,
where for $1\leq h\leq 3$, $l=3(s-1)+h$.
We here 
replace $P_t$ in (9) by
$$
P_t=\{(t,k,l,i,j)\ |\ 1\leq k\leq 3,\
l\in A^{*}_{tj},\ i\in B^{*}_{tjl},\ 1\leq j\leq n_{t}\}.
$$
Then $\eta(t,k,l,i,j)=x_{tk}z_{thq_s}f_{tp_i}z_{thq_s}g_{tj}$,
where $l=3(s-1)+h$ with $1\leq h\leq 3$ and $p_i\in B_{tj}(q_s)$.
We also replace $M_{t}$ by
$M_t=\{(l,i,j)\ |\ l\in A^{*}_{tj},\ i\in B^{*}_{tjl},\ 1\leq j\leq n_{t}\}$.
Let $\Sy=(V,E,E^{*})$ be as described in the proof of Lemma \ref{Mt},
where $V=M_t$ as above.
Then $\Sy$ is an R-colouring R-graph
with $\CG_{n_t}=\{V_1,\cdots,V_{n_t}\}$,
and the difference between this $\Sy$
and the one in the proof of Lemma \ref{Mt}
is simply that $V_j=\bigcup_{l=1}^{3}V_{jl}$
and $|V_{jl}|=m_t$ there
whereas $V_j=\bigcup_{l=1}^{3a_{tj}}V_{jl}$
and $|V_{jl}|=m_{tjs}$ here.
Since $a_{tj}>0$ and $m_{tjs}>0$,
we can easily see
that Theorem \ref{PS3} is also valid in this case
and that the same assertion as Lemma \ref{Mt} holds.
The remains of the proof are the same
as the proof of Theorem \ref{MainTh} (1).
$\Box$
\vskip12pt

If $F_{1}\subseteq F_{2}\subseteq\cdots$ are free groups,
then $F_{\infty}=\bigcup_{i=1}^{\infty}F_i$
contains a free subgroup $F$ with $|F|=|F_{\infty}|$.
In fact,
if either $|F_i|\leq\aleph_{0}$ for all $i$
or $|F_i|$ is a maximal cardinality for some $i$,
then the assertion is obvious.
Hence, it suffices to consider the case 
that their cardinalities are not bounded above.
We may then assume that
$|F_i|<|F_{i+1}|$ for all $i$.
Since each element of $F_i$ is a product of
finitely many basis elements of $F_{i+1}$,
each $F_{i+1}$ can be written as a free product
$G_{i+1}*H_{i+1}$,
where $G_{i+1}$ and $H_{i+1}$ are free subgroups of $F_{i+1}$
with $F_i\subseteq G_{i+1}$
and $|F_{i+1}|=|H_{i+1}|$.
Then $H_2*H_3*\cdots$ is a free subgroup of $F_{\infty}$
with the same cardinality as $F_{\infty}$.
Now,
it is well known
that a countable locally free group
is the union of an ascending sequence of free subgroups. 
Hence, by Theorem \ref{MainTh} (2),
we have

\begin{coro}\label{MainCor}
Let $F_1\subseteq F_2\subseteq\cdots \subseteq F_n\subseteq\cdots$
be an ascending chain of non-abelian free groups,
and $F_{\infty}=\cup_{i=1}^{\infty}F_i$.
Then the group ring $KF_{\infty}$ is primitive for any field $K$.
In particular,
every group ring of a countable non-abelian
locally free group over a field
is primitive.
\end{coro}

We are now in a position to
prove easily Theorem \ref{SubTh}:

\noindent
{\bf Proof.}\quad
[Proof of Theorem \ref{SubTh}]
By virtue of Lemma \ref{hnn},
we may assume that $\varphi(F)\ne F$.
Then $\Delta(F_{\varphi})=1$
by lemma \ref{fcc} (2).
Let $F_{i}$ be the subgroup of $F_{\varphi}$
generated by $\{ t^{i}ft^{-i}\ |\ f\in F\}$,
and $F_{\infty}=\bigcup_{i=1}^{\infty}F_{i}$.
By lemma \ref{fcc} (3),
$F_{\infty}$ is a normal subgroup of $F_{\varphi}$,
and it is also a locally free group
which is of type as described in Corollary \ref{MainCor}.
Hence, $KF_{\infty}$ is primitive
by Corollary \ref{MainCor}.
It is obvious that $F_{\varphi}/F_{\infty}$
is isomorphic to $\langle t\rangle$,
and thereby,
it follows from Lemma \ref{PBR} (1)
that $KF_{\varphi}$ is primitive.
$\Box$
\vskip12pt

Finally, we state the semiprimitivity
of group rings of ascending HNN extensions of free groups, 
which extends \cite[Corollary 3.7]{Ni07}
to the general cardinality case:

\begin{coro}\label{SubCor}
Let $F$ be a non-abelian free group,
and $F_{\varphi}$ the ascending HNN extension of
$F$ determined by $\varphi$.
If $K$ is any field then the group ring $KF_{\varphi}$ is semiprimitive.
\end{coro}

\noindent
{\bf Proof.}\quad
Let $K'$ be the prime field of $K$.
Since $|K'|\leq |F|$,
by virtue of Theorem \ref{SubTh},
$K'F_{\varphi}$ is primitive and so semiprimitive.
As is well known,
semiprimitive group rings are separable algebras,
thus semiprimitivity of group rings close under
extensions of coefficient fields,
and therefore $KF_{\varphi}$ is semiprimitive.
$\Box$
\vskip12pt


\end{document}